\documentclass[12pt]{article}
\usepackage{amsmath}
\usepackage{amssymb}
\usepackage{amsfonts}
\usepackage{amscd}

\def\no{\noindent}
\newtheorem{dfn}{Definition}[section]
\newtheorem{dfndesc}[dfn]{Definition-Description}
\newtheorem{rem}[dfn]{Remark}
\newtheorem{thm}[dfn]{Theorem}
\newtheorem{mthm}[dfn]{Main Theorem}
\newtheorem{lem}[dfn]{Lemma}
\newtheorem{sublem}[dfn]{Sublemma}
\newtheorem{prop}[dfn]{Proposition}
\newtheorem{propdfn}[dfn]{Proposition-Definition}
\newtheorem{cor}[dfn]{Corollary}

\newtheorem{ex}[dfn]{Example}

\newtheorem{cons}[dfn]{Consequence}

\newtheorem{metaques}[dfn]{Meta-Question}

\newtheorem{app}[dfn]{Application}
\newtheorem{ass}[dfn]{Assumption}
\newtheorem{add}[dfn]{Addendum}

\newtheorem{reform}[dfn]{Reformulation}
\newtheorem{refine}[dfn]{Refinement}

\def\proof{\par\medskip\noindent{\it Proof: }}

\def\qed{{\unskip\nobreak\hfil
        \penalty50\hskip1em\hbox{}\nobreak\hfil
        $\square$\parfillskip=0pt\finalhyphendemerits=0 \par\medskip}}
\def\restr{\mbox{\large \(|\)\normalsize}}
\def\R{\mathbb R}
\def\C{\mathbb C}

\def\H{\mathbb H}

\def\N{\mathbb N}

\def\al{\alpha}
\def\D{\partial}
\def\de{\delta}
\def\De{\Delta}

\def\eps{\epsilon}
\def\ga{\gamma}
\def\Ga{\Gamma}

\def\La{\Lambda}
\def\lra{\longrightarrow}

\def\Tau{{\cal T}}

\def\embed{\hookrightarrow}
\def\ol{\overline}
\def\om{\omega}

\def\ra{\rightarrow}

\def\si{\sigma}
\def\Si{\Sigma}
\def\tangle{\angle_{Tits}}
\def\tits{\partial_{Tits}}
\def\geo{\partial_{\infty}}
\def\karp{\partial_{\infty}^{fine}}
\def\karpD{\partial_{\infty}^{fine,\partial}}

\begin{document}
\title{A characterization of irreducible 
symmetric spaces and Euclidean buildings of higher rank 
by their asymptotic geometry}
\author{Bernhard Leeb}
\date{June 17, 1997}
\maketitle

\no
{\bf Abstract.}
We study geodesically complete and locally compact Hadamard spaces $X$ 
whose Tits boundary is a connected irreducible spherical building.
We show that $X$ is symmetric iff complete geodesics in $X$ do not branch
and a Euclidean building otherwise.
Furthermore,
every boundary equivalence
(cone topology homeomorphism preserving the Tits metric)
between two such spaces is induced by a homothety.
As an application,
we can extend the Mostow and Prasad rigidity theorems to
compact singular (orbi)spaces of nonpositive curvature
which are homotopy equivalent to a quotient of
a symmetric space or Euclidean building by a cocompact group of isometries.

\tableofcontents

\section{Introduction}

\subsection{Main result, background, motivation and an application}

Hadamard manifolds are simply-connected complete Riemannian 
manifolds of nonpositive sectional curvature. 
Prominent examples are Riemannian symmetric spaces of noncompact type,
but many more examples occur as universal covers of closed nonpositively
curved manifolds.
For instance, most Haken 3-manifolds admit metrics of nonpositive curvature
\cite{thesis}. 
Not the notion of sectional curvature itself, however the notion of 
an upper curvature bound can be expressed purely by inequalities 
involving the distances between finitely many points
but no derivatives of the Riemannian metric, 
and hence generalizes from the narrow world of Riemannian manifolds 
to a wide class of metric spaces, cf.\ \cite{Aleks}. 
The natural generalization of Hadamard manifolds 
are Hadamard spaces,
i.e.\ complete geodesic metric spaces which are 
nonpositively curved in the (global) sense of distance comparison, 
see \cite{Ballmann,symm}. 
Hadamard spaces comprise besides Hadamard manifolds a large class 
of interesting singular spaces, among them 
Euclidean buildings (the discrete cousins of symmetric spaces), 
many piecewise Euclidean or Riemannian complexes occuring,
for instance, in geometric group theory, 
and branched covers of Hadamard manifolds. 
Hadamard spaces received much attention in the last decade, 
notably with view to geometric group theory, 
a main impetus coming from Gromov's work 
\cite{HypGrp,AsyInv}. 

We recall that a fundamental feature of a Hadamard space is 
the convexity of its distance function 
with the drastic consequences such as uniqueness of geodesics 
and in particular contractibility. 
This illustates that already geodesics, 
undoubtedly fundamental objects in geometric considerations, 
are rather well-behaved and their behavior 
can be to some extent controlled, 
which gets the foot in the door for a more advanced geometric understanding. 
The importance of the geometry of nonpositive curvature 
lies in the coincidence that one has a rich supply of interesting examples
reaching into many different branches of mathematics 
(like geometric group theory, representation theory, arithmetic)
and, at the same time, these spaces share simple basic geometric properties 
which makes them understandable to a certain extent and in a uniform way. 

\medskip 
We will be interested in asymptotic information and the restrictions which 
it imposes on the geometry of a Hadamard space $X$. 
This is related to the rigidity question, 
already classical in global Riemannian geometry, 
how topological properties of a (for instance closed) Riemannian manifold 
with certain local (curvature) constraints 
are reflected in its geometry\footnote{
Since the universal cover is contractible, the entire topological information 
is contained in the fundamental group and one can ask which of its 
algebraic properties are visible in the geometry.} 
and below (\ref{ZusatzzuHauptsatz}) 
we will present an application in this direction. 

Let us first describe which asymptotic information we consider. 
The {\em geometric} or {\em ideal boundary} $\geo X$ 
of a Hadamard space $X$ is defined as 
the set of equivalence classes of asymptotic geodesic rays.\footnote{
For two unit speed geodesic rays $\rho_1,\rho_2:[0,\infty)$ 
the distance $d(\rho_1(t),\rho_2(t))$ of travellers along the rays
is a convex function. 
If it is bounded (and hence non-increasing) 
the rays are called {\em asymptotic}.}
The topology on $X$ extends to a natural {\em cone topology} 
on the geometric completion $\bar X=X\cup\geo X$ which is compact 
iff $X$ is locally compact. 
The ideal boundary points $\xi\in\geo X$ can be thought of as the ways
to go straight to infinity\footnote{
Examples: The geometric completion of hyperbolic plane can be obtained 
by taking the closure in the Poincar\'e disk model;
one obtains the Poincar\'e Compact Disk model. 
The geometric boundary of a metric tree is the set of its ends
which is a Cantor set if it has no isolated points.}.
It is fair to say that the topological type of $\geo X$ is not 
a very strong invariant, for example it is a $(n-1)$-sphere for 
any $n$-dimensional Hadamard manifold. 

Besides the cone topology there is another interesting structure 
on $\geo X$, 
namely the {\em Tits angle metric} introduced by Gromov 
in full generality in \cite{BGS}. 
For two points $\xi_1,\xi_2\in\geo X$ at infinity
their Tits angle $\tangle(\xi_1,\xi_2)$ measures the 
maximal visual angle $\angle_x(\xi_1,\xi_2)$ under which they can be seen
from a point $x$ inside $X$, 
or equivalently,
it measures the asymptotic linear rate at which unit speed geodesic rays
$\rho_i$ asymptotic to the ideal points $\xi_i$ diverge from each other. 
If $X$ has a strictly negative curvature bound
the {\em Tits boundary} $\tits X=(\geo X,\tangle)$ 
is a discrete metric space and only of modest interest. 
However, if $X$ features substructures of extremal curvature zero,
such as flats, i.e.\ convex subsets isometric to Euclidean space, 
then connected components appear in the Tits boundary 
and the Tits metric becomes an interesting structure.\footnote{
$\tits SL(3,\R)/SO(3)$ is the 1-dimensional spherical building 
associated to the real projective plane.} 
The cone topology together with the Tits metric on $\geo X$ are the 
{\em asymptotic data} which we consider here. 
Our results find shelter under the roof of the following: 

\begin{metaques}
What are the implications of these asymptotic data for the geometry 
of a Hadamard space? 
\end{metaques}

The main result is the following characterization of symmetric spaces and 
Euclidean buildings of higher rank as Hadamard spaces with 
spherical building boundary: 

\begin{mthm}
\label{hauptsatz}
Let $X$ be a locally compact Hadamard space 
with extendible geodesic segments\footnote{I.e.\ every geodesic segment 
is contained in a complete geodesic.} 
and assume that $\tits X$ is a connected thick irreducible spherical building.
Then $X$ is a Riemannian symmetric space or a Euclidean building.
\end{mthm}

In the smooth case, i.e.\ for Hadamard manifolds, 
\ref{hauptsatz} follows from work of Ballmann and Eberlein,
cf.\ \cite[Theorem B]{symmdiff}, or else from arguments of 
Gromov \cite{BGS} and Burns-Spatzier \cite{BurnsSpatzier}. 
There is a dichotomy into two cases,
according to whether geodesics in $X$ branch or not. 
In the absence of branching the ideal boundaries are very symmetric 
because there is an involution $\iota_x$ of $\geo X$ 
at every point $x\in X$, 
and one can adapt arguments from Gromov in the proof of his rigidity theorem 
\cite{BGS}. 
Our main contribution lies in the case of geodesic branching. 
There the boundary at infinity admits in general no non-trivial symmetries 
and another approach is needed. 

We show moreover that for the spaces considered in \ref{hauptsatz} 
the extreme situation occurs that 
$X$ is completely determined by its asymptotic data up to a scale factor: 

\begin{add}
\label{ZusatzzuHauptsatz}
Let $X$ be a 
symmetric space or a thick Euclidean building, 
irreducible and of rank $\geq2$,
and let $X'$ be another such space. 
Then any boundary isomorphism
(cone topology homeomorphism preserving the Tits metric)
\begin{equation}
\label{RandIso}
\phi:\geo X\ra\geo X'
\end{equation}
is induced by a homothety.
\end{add}

\ref{ZusatzzuHauptsatz} 
follows from Tits classification for automorphisms of spherical buildings 
in the cases when $X$ has many symmetries, 
e.g.\ 
when it is a Riemannian symmetric space or a Euclidean building associated 
to a simple algebraic group over 
a local field with non-archimedean valuation. 
This is in particular true if $rank(X)\geq3$ 
however it does not cover the cases when $X$ is a rank 2 Euclidean building 
with small isometry group. 
Our methods provide a uniform proof in all cases and in particular 
a direct argument in the symmetric cases.

\bigskip
A main motivation for us was Mostow's Strong Rigidity Theorem 
for locally symmetric spaces,
namely the irreducible case of higher rank: 

\begin{thm}[\cite{Mostow}]
Let $M$ and $M'$ be locally symmetric spaces whose universal covers 
are irreducible symmetric spaces of rank $\geq2$. 
Then any isomorphism $\pi_1(M)\ra\pi_1(M')$ of fundamental groups 
is induced by a homothety $M\ra M'$. 
\end{thm}

It is natural to ask whether locally symmetric spaces are rigid 
in the wider class of closed manifolds of nonpositive sectional 
curvature. 
This is true and the content of Gromov's Rigidity Theorem \cite{BGS}. 
As an application of our main results 
we present an extension of Mostow's theorem 
as well as Prasad's analogue for compact quotients of Euclidean buildings 
\cite{Prasad} 
to the larger class of singular nonpositively curved (orbi)spaces:

\begin{app}
\label{verallgMostow}
Let $X$ be a locally compact Hadamard space with extendible geodesic 
segments and
let $X_{model}$ be a symmetric space (of noncompact type) or a thick Euclidean
building.
Suppose furthermore that all irreducible factors of $X_{model}$ have rank
$\geq2$.
If the same finitely generated group $\Ga$ acts
cocompactly and properly discontinuously on $X$ and $X_{model}$
then, after suitably rescaling the metrics on the irreducible factors of
$X_{model}$, there is a $\Ga$-equivariant isometry
$X\ra X_{model}$.
\end{app}

I.e.\ among (possibly singular) geodesically complete 
compact spaces of nonpositive curvature (in the local sense), 
quotients of irreducible higher rank symmetric spaces or Eulidean buildings 
are determined by their homotopy type. 

\begin{ex}
On a locally symmetric space with irreducible higher rank universal cover
there exists no piecewise Euclidean singular metric of nonpositive curvature. 
\end{ex} 

As we said, 
\ref{verallgMostow} is due to Gromov \cite{BGS} 
in the case that $X$ is smooth Riemannian. 
Although we extend Gromov's extension of Mostow Rigidity 
further to singular spaces, 
the news of \ref{verallgMostow} lie mainly in the building case.

\bigskip
{\bf Acknowledgements.}
I would like to use this opportunity to thank 
the Sonderforschungsbereich SFB 256 and Werner Ballmann 
for generously supporting my mathematical research over the past years. 

I make extensive use of the geometric approach 
to the theory of generalized Tits buildings developed in \cite{symm}
together with Bruce Kleiner 
and thank him for numerous discussions  
related to the subject of this paper. 
I am grateful to Werner Ballmann and Misha Kapovich 
for pointing out countably many errors in a late version of the manuscript, 
to Jens Heber for communicating a reference 
and to the Complutense library in Madrid for several crucial suggestions. 
Part of this work was done while visiting the 
Instituto de Matematicas de la UNAM, 
Unidad Cuernavaca, in Mexico.
I thank the institute for its great hospitality,
especially Pepe Seade and Alberto Verjovski.

\subsection{Around the argument}

In this section I attempt to describe the scenery around the proofs of 
\ref{hauptsatz} and \ref{ZusatzzuHauptsatz}  
for the rank 2 case,
i.e.\  when maximal flats in $X$
have dimension 2 and $\tits X$ is 1-dimensional. 
The rank 2 case is anyway the critical case 
because there the rigidity is qualitatively weaker 
than in the case of rank $\geq3$.
This difference is reflected in 
Tits classification theorem for spherical buildings
\cite{Tits} which asserts, roughly speaking, 
that all thick irreducible spherical buildings of 
rank $\geq 3$ (that is, dimension $\geq2$) 
are canonically attached to simple algebraic or classical groups. 
In contrast there exist uncountably many absolutely asymmetric 1-dimensional 
spherical buildings, 
for example those corresponding to exotic projective planes, 
and uncountably many of them occur as Tits boundaries 
of rank 2 Euclidean buildings with trivial isometry group. 

\medskip
For a singular geodesic $l$ in $X$, 
i.e.\ a geodesic asymptotic to vertices in $\tits X$, 
we consider the union $P(l)$ of all geodesics parallel to $l$ 
and its cross section $CS(l)$ which is a locally compact Hadamard space 
with discrete Tits boundary 
(as for a rank 1 space).\footnote{
For instance, if $X=SL(3,\R)/SO(3)$ 
then the cross sections of singular geodesics are hyperbolic planes.
More generally, if $X$ is a symmetric space of rank 2 then these 
cross sections are rank-1 symmetric spaces.
If $X$ is a Euclidean building of rank 2 they are 
rank-1 Euclidean buildings, i.e.\ metric trees.}
For two asymptotic geodesics $l$ and $l'$ 
one can canonically identify the ideal boundaries $\geo CS(l)$ 
and $\geo CS(l')$. 
A priori this identification of boundaries does not extend to an isometry 
between the cross sections, 
but it extends to an isometry between their ``convex cores'',
i.e.\ between closed convex subsets 
$C\subseteq CS(l)$ and $C'\subseteq CS(l')$ 
which are minimal among the closed convex subsets satisfying 
$\geo C=\geo CS(l)$ and $\geo C'=\geo CS(l')$. 
Due to a basic rigidity phenomenon in the geometry of nonpositively curved 
spaces, 
the so-called ``Flat Strip Theorem'', 
the convex cores are unique up to isometry. 
We see that,
essentially due to the connectedness of the Tits boundary, 
there are 
many natural identifications between the various parallel sets
and observe that, by composing them, 
one can generate large groups of isometries acting 
on the cores $C_l$ of the cross sections 
(see section \ref{holonomy}). 
We denote the closures of these ``holonomy'' subgroups in $Isom(C_l)$ by 
$Hol(l)$. 
They are large in the sense that 
$Hol(l)$ acts 2-fold transitively on $\geo CS(l)$.\footnote{
If $X$ is a rank-2 symmetric space, $Hol(l)$ contains 
the identity component of the isometry group of the rank-1 symmetric 
space $CS(l)$. 
So in the example $X=SL(3,\R)/SO(3)$ ($X=SL(3,\C)/SU(3)$) 
the action of $Hol(l)$ on the boundary of the hyperbolic plane 
(hyperbolic 3-space) $CS(l)$ is even 3-fold transitive 
(by M\"obius transformations). 
More generally the action is 3-fold transitive if $\tits X$
is the spherical building associated to an (abstract) projective plane.} 
Hence, 
{\em however unsymmetric $X$ itself may be, the cross sections 
of its parallel sets are always highly symmetric,} 
and this is the key observation at the starting point of our argument. 

The high symmetry imposes a substantial restriction 
on the geometry of the cross sections 
and the major step in our proof of \ref{hauptsatz} is a rank-1 analogue for 
spaces with high symmetry: 

\begin{thm}
\label{Rang1Analog}
Let $Y$ be a locally compact Hadamard space with extendible rays 
and at least 3 points at infinity. 
Assume that $Y$ contains a closed convex subset $C$ with full ideal boundary
$\geo C=\geo Y$ 
so that $Isom(C)$ acts 2-fold transitively on $\geo C$. 
Then the following dichotomy occurs:

1.\ 
If some complete geodesics in $Y$ branch 
then $Y$ is isometric to the product of a metric tree 
(with edges of equal length) 
and a compact Hadamard space. 

2.\ 
If complete geodesics in $Y$ do not branch then 
there exists a rank-1 Riemannian symmetric model space $Y_{model}$,
and a boundary homeomorphism $\geo Y\ra\geo Y_{model}$ 
carrying $Isom_o(C)$ to $Isom_o(Y_{model})$.\footnote{
It seems unclear whether in this case one should not be able to find 
an embedded rank-1 symmetric space inside $Y$.}
\end{thm}

In particular, 
the ideal boundary of every cross section is homeomorphic 
to a sphere, a Cantor set or a finite set of cardinality $\geq3$. 

\bigskip
As we explained, the geometry of $X$ 
is rigidified by the various identifications between cores of cross sections 
of parallel sets. 
This can be nicely built in the picture of the geometric
compactification of $X$ as follows: 
We mentioned that the convex cores of the cross sections $CS(l)$ 
for all lines $l$ asymptotic to the same vertex $\xi\in\tits X$ 
can be canonically identified to a Hadamard space $C_{\xi}$. 
It has rank 1 in the sense that it satisfies the visibility property, 
or equivalently, its Tits boundary is discrete. 
$\geo C_{\xi}$ can be reinterpreted as the compact topological space 
of Weyl chambers (arcs) emanating from the vertex $\xi$. 
One can now blow up the geometric boundary $\geo X$ 
by replacing each vertex $\xi$ by the geometric compactification 
$\bar C_{\xi}$ and gluing the endpoints of Weyl arcs to the corresponding
boundary points in $\geo C_{\xi}$. 
This generalizes a construction by Karpelevi\v{c} for symmetric spaces
\cite{karpel}. 
We denote the resulting refined boundary by $\karp X$, 
and by $\karpD X$ 
the part which one obtains by inserting only the boundaries
$\geo C_{\xi}$ instead of the full compactifications $\bar C_{\xi}$. 
{\em The rigidity expresses itself in the action of the holonomy groupoid 
which appears on the blown up locus of the refined boundary $\karp X$
due to the connectedness of $\tits X$:} 
For any two antipodal vertices $\xi_1,\xi_2\in\tits X$, 
i.e.\ vertices of Tits distance $\tangle(\xi_1,\xi_2)=\pi$, 
there is a canonical isometry 
\begin{equation}
\label{IdentQuersch}
C_{\xi_1}\leftrightarrow C_{\xi_2}
\end{equation}
because the spaces $C_{\xi_i}$ embed as minimal convex subsets 
into the cross section $CS(\{\xi_1,\xi_2\})$ of the family of 
parallel geodesics asymptotic to $\xi_1,\xi_2$.
We can compose such isometries hopping along finite sequences 
of successive antipodes.
For any two vertices $\xi,\eta\in\tits X$ 
we denote by $Hol(\xi,\eta)\subseteq Isom(C_{\xi},C_{\eta})$ 
the closure of the subset of all isometries $C_{\xi}\ra C_{\eta}$ 
which arise as finite composites of isometries (\ref{IdentQuersch})
(cf.\ section \ref{holonomy}). 
In particular, 
the holonomy groups $Hol(\xi):=Hol(\xi,\xi)$ act on the inserted 
spaces $\bar C_{\xi}$. 
These actions can be thought of as an additional geometric structure 
on the spaces $\geo C_{\xi }$, 
namely as the analogue of a conformal structure; 
for instance if $\geo C_{\xi}$ is homeomorphic to a sphere 
then due to \ref{Rang1Analog} it can be identified with the boundary 
of a rank-1 symmetric space up to conformal diffeomorphism. 

\medskip
{\em Comment on the proof of \ref{hauptsatz}:} 
It is easy to see 
that all cross sections $CS(l)$ have extendible geodesic rays
if $X$ is geodesically complete (\ref{QuerschnittehabenerweiterbareStrahlen}). 

{\em If some complete geodesics branch} in $X$ 
then there is a cross section $CS(l)$ with branching geodesics 
(\ref{dannverzweigenauchdiequerschnitte}) and, 
apparently less trivially to verify, 
even all cross sections have this property (\ref{einerdannallecantor}). 
The rank-1 result \ref{Rang1Analog} then implies that the cross sections 
of all parallel sets are metric trees
(up to a compact factor). 
From this point it is fairly straight-forward 
to conclude in one way or another that $X$ is a Euclidean building 
(section \ref{BeweisvonHauptsatz}). 

{\em If complete geodesics in $X$ do not branch}
we can adapt arguments of Gromov 
from the proof of his Rigidity Theorem \cite{BGS}. 
The reflections at points $x\in X$ give rise to involutive automorphisms 
$\iota_x:\geo X\ra\geo X$ of the topological spherical building $\geo X$. 
One obtains a proper map $X\embed Aut(\geo X)$ into the group of
boundary automorphisms and hence finds oneself in the situation that 
the topological spherical building $\tits X$ is highly symmetric.
(It satisfies the so-called Moufang property.) 
$Aut(\geo X)$ is a locally compact topological group \cite{BurnsSpatzier}.
Similar to \cite{BurnsSpatzier}, 
establishing transitivity and contraction properties for the dynamics 
of $Aut(\geo X)$ on $\geo X$ allows to show, 
using a deep result by Gleason and Yamabe on the approximation 
of locally compact topological groups by Lie groups, 
that $Aut(\geo X)$ is a semisimple Lie group and the isometry group 
of a Riemannian symmetric model space $X_{model}$. 
The involutions $\iota_x$ can be characterized as order 2 elements
with compact centralizer and hence correspond to point reflections 
in $X_{model}$. 
One obtains a map $\Phi:X\ra X_{model}$ which is clearly affine in the sense
that it preserves flats. 
It immediately follows that $\Phi$ is a homothety,
concluding the proof of \ref{hauptsatz}. 

\medskip
{\em Comment on the proof of \ref{ZusatzzuHauptsatz}:}
Any boundary isomorphism (\ref{RandIso})
has continuous differentials
\begin{equation}
\label{RandDiff}
\Si_{\xi}\phi:\Si_{\xi}\tits X\ra\Si_{\phi\xi}\tits X'
\end{equation}
and hence lifts to a map 
\begin{equation}
\label{TeilLift}
\karpD X\ra\karpD X'
\end{equation}
of partially refined boundaries.
The differentials (\ref{RandDiff})
are conformal in the sense that they preserve the holonomy action,
i.e.\ the induced homeomorphisms
\[ Homeo(\geo C_{\xi},\geo C_{\eta})\ra
Homeo(\geo C_{\phi\xi},\geo C_{\phi\eta}) \]
carry the holonomy groupoid $Hol^X$ to the holonomy groupoid $Hol^{X'}$.
This sets us on the track towards the proof of \ref{ZusatzzuHauptsatz}:
After proving \ref{hauptsatz}
we may assume that $X$ is a symmetric space or a Euclidean building.
Then the $C_{\xi}$ are rank-1 symmetric spaces or metric trees,
respectively, 
and the differentials $\Si_{\xi}\phi$ actually extend to homotheties 
$C_{\xi}\ra C_{\phi\xi}$.
This means that the lift 
(\ref{TeilLift}) of $\phi$ improves to a holonomy equivariant map 
\begin{equation}
\label{induziertKarp}
\karp X\lra\karp X' 
\end{equation}
between the full refined geometric boundaries. 
Since points in the blow ups $C_{\xi}$ are equivalence classes 
of strongly asymptotic geodesics, 
(\ref{induziertKarp}) encodes a correspondence between singular 
geodesics in $X$ and $X'$. 
If $X$ is a Euclidean building then this can be used 
in a final step to set up a correspondence between vertices 
which preserves apartments and extends to a homothety $X\ra X'$,
hence concluding the proof of \ref{ZusatzzuHauptsatz} in this case
(section \ref{fliege}) .
If $X$ is a symmetric space then \ref{ZusatzzuHauptsatz} 
already follows from the arguments in the proof of \ref{hauptsatz}. 

\medskip
The paper is desorganized as follows: 
In section \ref{prelim} we discuss preliminaries.
In particular we establish the existence of convex cores 
for Hadamard spaces under fairly general conditions 
(section \ref{secconvexcores}) and introduce the
spaces of strong asymptote classes which will serve as an important 
tool in the construction of the holonomy groupoid. 
The holonomy groupoid is discussed 
in section \ref{holonomy} 
where we explain the symmetries of parallel sets.
In section \ref{secrangeins} 
we prove the rigidity results for ``rank 1'' spaces with high symmetry
and 
in section \ref{sectionhadamardspaceswithbuildingboundary} 
the main results for higher rank spaces.

\section{Preliminaries}
\label{prelim}

\subsection{Hadamard spaces}

For basics on Hadamard spaces and, more generally, 
spaces with curvature bounded above we refer to 
the first two chapters of \cite{Ballmann} and 
section 2 of \cite{symm}. 
Spaces of directions and Tits boundaries are discussed there 
and it is verified that they are CAT(1) spaces. 
Let us emphasize that we mean by the {\em Tits boundary} $\tits X$ 
of the Hadamard space $X$ the geometric boundary $\geo X$ 
equipped with the {\em Tits angle metric} $\tangle$
and not with the associated path metric\footnote{
If $\tits X$ is a spherical building then it has diameter $\pi$ 
with respect to the path metric and hence the path metric coincides 
with $\tangle$.}. 
 
In the following paragraphs 
we supply a few auxiliary facts needed later in the text. 

\subsubsection{Filling spheres at infinity by flats}

The following result generalizes an observation by Schroeder 
in the smooth case, cf.\ \cite{BGS}. 

\begin{prop}
\label{sphaerenfuellen}
Let $X$ be a locally compact Hadamard space 
and let $s\subseteq\tits X$ be a unit sphere which does not bound 
a unit hemisphere in $\tits X$. 
Then there exists a flat $F\subseteq X$ with $\geo F=s$.
\end{prop}
\proof
Let $s$ be isometric to the unit sphere of dimension $d\geq0$ 
and pick $d+1$ pairs of antipodes $\xi_0^{\pm},\dots,\xi_d^{\pm}$ 
so that 
\begin{equation}
\label{winkelpihalbe}
\tangle(\xi_i^{\pm},\xi_j^{\pm})=\pi/2
\qquad\hbox{ and }\qquad
\tangle(\xi_i^{\pm},\xi_j^{\mp})=\pi/2
\end{equation} 
for all $i\neq j$. 
If for some point $x\in X$ and some index $i$ holds 
$\angle_x(\xi_i^+,\xi_i^-)=\pi$ then the union $X'=P(\{\xi_i^+,\xi_i^-\})$ 
of geodesics asymptotic to $\xi_i^{\pm}$ is non-empty 
and $s$ determines a $(d-1)$-sphere $s'\subseteq\tits X'$ which does not 
bound a unit hemisphere. 
Moreover any flat $F'\subseteq X'$ filling $s'$ determines a flat $F$
filling $s$ and we are reduced to the same question 
with one dimension less. 
We can hence proceed by induction on the dimension $d$ 
and the claim follows if we can rule out the situation that 
\begin{equation}
\label{winkelkleinerpihalbe}
\angle_x(\xi_i^+,\xi_i^-)<\pi
\end{equation}
holds for all $x$ and $i$. 
In this case we obtain a contradiction as follows. 
Assume first that for some (and hence any) point $x_0\in X$ the 
intersection of the horoballs $Hb(\xi_i^{\pm},x_0)$ is unbounded
and thus contains a complete geodesic ray $r$.
The ideal endpoint $\eta\in\geo X$ of $r$ satisfies 
$\tangle(\eta,\xi_i^{\pm})\leq\pi/2$ 
because the Busemann functions 
$B_{\xi_i^{\pm}}$ monotonically non-increase along $r$. 
By the triangle inequality follows 
$\tangle(\eta,\xi_i^{\pm})=\pi/2$ 
because $\xi_i^{\pm}$ are antipodes. 
The CAT(1) property of $\tits X$ then implies that there is a unit hemisphere 
$h\subseteq\tits X$ with center $\eta$ and boundary $s$,
but this contradicts our assumption. 
Therefore the intersection of the horoballs $Hb(\xi_i^{\pm},x_0)$ 
is compact for all $x_0\in X$
and the convex function $\max B_{\xi_i^{\pm}}$ is proper 
and assumes a minimum in some point $x$. 
Denote by $r_i^{\pm}:[0,\infty)\ra X$ 
the ray with $r_i^{\pm}(0)=x$ and $r_i^{\pm}(\infty)=\xi_i^{\pm}$. 
(\ref{winkelpihalbe}) implies that $B_{\xi_i^{\pm}}$ non-increases
along $r_j^{\pm}$ for $i\neq j$. 
Hence, if
$x_j$ denotes the midpoint of the segment $\ol{r_j^+(1)r_j^-(1)}$ 
then $B_{\xi_i^{\pm}}(x_j)\leq B_{\xi_i^{\pm}}(x)$ for all $i$ 
and, by (\ref{winkelkleinerpihalbe}), 
$B_{\xi_j^{\pm}}(x'_j) < B_{\xi_j^{\pm}}(x)$ 
for some point $x'_j\in\ol{xx_j}$. 
This means that by replacing $x$ we can decrease the values of one pair of
Busemann functions while not increasing the others.
By iterating this procedure at most $d+1$
times we find a point $x'$ with 
$\max B_{\xi_i^{\pm}}(x') < \max B_{\xi_i^{\pm}}(x)$,
a contradiction. 
\qed

\subsubsection{Convex cores}
\label{secconvexcores}

For a subset $A\subseteq\geo Y$ we denote by ${\cal C}_A$ the family of
closed convex subsets $C\subseteq Y$ with $\geo C\supseteq A$.
${\cal C}_A$ is non-empty, partially ordered and closed under intersections.

\begin{prop}
\label{exuniqueminimalconvexsubset}
Let $Y$ be a locally compact Hadamard space.

1.\  
Suppose that $s\subseteq A\subseteq\geo Y$ 
and $s$ is a unit sphere with respect to the Tits metric 
which does not bound a unit hemisphere. 
Then ${\cal C}_A$ contains a minimal element.

2.\ 
Suppose that $A\subseteq\geo Y$ so that ${\cal C}_A$ 
has minimal elements. 
Then the union $Y_0$ of all minimal elements in ${\cal C}_A$ is a 
convex subset of $Y$.
It decomposes as a metric product
\begin{equation}
\label{flatstripofconvexsets}
Y_0\cong C\times Z
\end{equation}
where $Z$ is a compact Hadamard space
and the layers $C\times\{z\}$ are precisely the minimal elements
in ${\cal C}_A$.
\end{prop}
\proof
According to \ref{sphaerenfuellen}, 
there exists a non-empty family of flats in $Y$ 
with ideal boundary $s$ 
and the family is compact because otherwise $s$ would bound 
a unit hemisphere. 
The union $P(s)$ of these flats is a convex subset of $Y$. 

\begin{sublem}
Let $F$ be a flat and $C$ a closed convex subset in $Y$ so that 
$\geo F\subseteq\geo C$.
Then $C$ contains a flat $F'$ parallel to $F$.
\end{sublem}
\proof
For any points $x\in C$ and $y\in F$ there is a point $x'\in C$ 
so that $d(x',y)\leq d(x,F)$. 
Hence there exists a point $x''\in C$ which realizes the nearest point
distance of $F$ and $C$:
$d(x'',F)=d(C,F)$. 
Then the union of rays emanating from $x''$ and asymptotic to points 
in $\geo F$ 
forms a flat $F'$ parallel to $F$.
\qed

Hence every convex subset $C\in {\cal C}_A$
intersects $P(s)$ in a non-empty compact family of flats 
and therefore determines a non-empty compact subset $U(C)$ 
in the compact cross section $CS(s)$ (compare definition \ref{DefParallelset}).
We order the sets $C\in {\cal C}_A$ by inclusion and observe that 
the assignment $C\mapsto U(C)$ preserves inclusion. 

\begin{sublem}
Let $(S_{\iota})$ be an ordered decreasing family of non-empty compact 
subsets of a compact metric space $Z$. 
Then the intersection of the $S_{\iota}$ is not empty.
\end{sublem}
\proof
For every $n\in\N$ we can cover $Z$ by finitely many balls of radius $1/n$ 
and therefore 
there exists a ball $B_{1/n}(z_n)$ which intersects 
all sets $S_{\iota}$. 
Any accumulation point of the sequence $(z_n)$ is contained 
in the intersection of the $S_{\iota}$.
\qed

Any decreasing chain of sets $C_{\iota}\in {\cal C}_A$ yields 
a decreasing chain of compact cross sections $U(C_{\iota})$
and hence has non-empty intersection. 
It follows that $\emptyset\neq\bigcap C_{\iota}\in {\cal C}_A$ and,
by Zorn's lemma or otherwise,  
we conclude that ${\cal C}_A$ contains a minimal non-empty subset. 

Now let $C_1,C_2\in{\cal C}_A$ be minimal.
For any $y_1\in C_1$ the closed convex subset
$\{y\in C_1:d(y,C_2)\leq d(y_1,C_2)\}$ of $C_1$ contains $A$ in its 
ideal boundary
and, by minimality of $C_1$, is all of $C_1$.
It follows that $d(\cdot,C_2)$ is constant on $C_1$
and the nearest point projection $p_{C_2C_1}:C_1\to C_2$ is an isometry.
For a decomposition $d(C_1,C_2)=d_1+d_2$ as a sum of positive numbers,
the set
$\{y\in Y:d(y,C_i)=d_i\hbox{ for $i=1,2$}\}$ is a minimal element in
${\cal C}_A$.
Hence $Y_0$ is convex.

\begin{sublem}
For minimal elements $C_1,C_2,C_3\in{\cal C}_A$ the self-isometry
$\psi=p_{C_1C_2}\circ p_{C_2C_3}\circ p_{C_3C_1}$ of $C_1$ is the identity.
\end{sublem}
\proof
$\psi$ preserves the central flat $f$ in $C_1$ with ideal boundary $s$.
Furthermore, 
$\psi\restr_f$ 
preserves all Busemann functions centered at ideal points $\in s$. 
Thus $\psi$ restricts to the identity on $f$. 
Since $\geo\psi=id$ and $C_1$ is minimal it follows that 
$\psi$ fixes $C_1$ pointwise. 
\qed

Now choose a minimal set $C\in{\cal C}_A$ and a point $y\in C$.
Then the set $Z$ of points $p_{C'C}(y)$, 
where $C'$ runs through all minimal elements in ${\cal C}_A$, 
is convex.
It is easy to see that $Y_0$ is canonically isometric to $C\times Z$.
$Z$ must be compact because $CS(s)$ is. 
This concludes the proof of \ref{exuniqueminimalconvexsubset}. 
\qed

The compact Hadamard space $Z$
in (\ref{flatstripofconvexsets})
has a well-defined center $z_0$.
We call the layer $C\times\{z_0\}$ the {\em central} minimal convex subset
in ${\cal C}_A$.

\begin{dfn}
If ${\cal C}_{\geo Y}$ has minimal elements then 
the {\bf convex core} $core(Y)$ of $Y$ is defined as the central 
minimal closed convex subset in 
${\cal C}_{\geo Y}$.
\end{dfn}

If the convex core exists it is preserved by all isometries of $Y$.

\begin{lem}
\label{nontranslatable}
Let $Y$ be a locally compact Hadamard space
which has a convex core. 
If $core(Y)$ has no Euclidean factor 
then any isometry with trivial action at infinity
fixes $core(Y)$ pointwise.
\end{lem}
\proof
Let $\phi$ be an isometry which acts trivially at infinity.
Then its displacement function is constant on the central convex subset $C$.
It is zero because $C$ does not split off a Euclidean factor.
\qed

\subsubsection{Spaces of strong asymptote classes} 
\label{secstrongasy}

Let $X$ be a Hadamard space.
For a point $\xi\in\geo X$ let us consider the rays asymptotic to $\xi$.
The {\em asymptotic distance} of 
two rays $\rho_i:[0,\infty)\ra X$ is given by 
their nearest point distance 
\begin{equation}
\label{asymptotischeDistanz} 
d_{\xi}(\rho_1,\rho_2)=\inf_{t_1,t_2\ra\infty}
d(\rho_1(t_1),\rho_2(t_2)), 
\end{equation}
which equals $\lim_{t\ra\infty}d(\rho_1(t),\rho_2(t))$
when the rays are parametrized so that 
$B_{\xi}\circ\rho_1\equiv B_{\xi}\circ\rho_2$. 
We call the rays $\rho_i$ {\em strongly asymptotic} 
if their asymptotic distance is zero. 
The asymptotic distance (\ref{asymptotischeDistanz}) 
defines a metric on the space $X^*_{\xi}$ of strong asymptote classes. 

\begin{propdfn}
\label{DefstarkeAsyklassen}
The metric completion $X_{\xi}$ of $X^*_{\xi}$ is a Hadamard space.
\end{propdfn}
\proof
Any two points in $X^*_{\xi}$ are represented by rays 
$\rho_1,\rho_2:[0,\infty)\ra X$ asymptotic to $\xi$ and initiating 
on the same horosphere centered at $\xi$. 
Denote by $\mu_s:[s,\infty)\ra X$ the ray asymptotic to $\xi$ whose
starting point $\mu_s(s)$ is the midpoint of $\ol{\rho_1(s)\rho_2(s)}$. 
The triangle inequality implies that 
$d(\rho_1(t),\mu_s(t))+d(\mu_s(t),\rho_2(t))-d(\rho_1(t),\rho_2(t))
\leq d(\rho_1(s),\rho_2(s)) -d(\rho_1(t),\rho_2(t))
\ra 0$
as $s,t\ra\infty$ with $s\leq t$. 
Hence $d(\mu_s(t),\mu_t(t))\to0$ and $d_{\xi}(\mu_s,\mu_t)\to0$,
i.e.\ $(\mu_s)$ is a Cauchy sequence and its limit in $X_{\xi}$
is a midpoint for $[\rho_1]$ and $[\rho_2]$. 
In this manner we can assign to every pair of points 
$[\rho_1],[\rho_2]\in X_{\xi}^*$ 
a well-defined midpoint $m\in X_{\xi}$.
If $[\rho'_1],[\rho'_2]\in X_{\xi}^*$ is another pair of points 
so that $d([\rho_i],[\rho'_i])\leq\de$ then $d(m,m')\leq\de$. 
It follows that there exist midpoints for all pairs of points in $X_{\xi}$.
As a consequence, 
any two points in $X_{\xi}$ can be connected by a geodesic. 

Any finite configuration ${\cal F}$ 
of points in $X^*_{\xi}$ corresponds to
a finite set of rays $\rho_i:[0,\infty)\ra X$ asymptotic to $\xi$
and synchronized so that for any time $t$ the set ${\cal F}_t$ 
of points $\rho_i(t)$ lies on one horosphere centered at $\xi$. 
The finite metric spaces $({\cal F}_t,d_X)$ Hausdorff converge 
to $({\cal F},d_{\xi})$ 
and hence distance comparison inequalities are inherited. 
It follows that geodesic triangles satisfy the 
CAT(0) comparison inequality. 
\qed

We will also $X_{\xi}$ call the {\em space of strong asymptote classes}
at $\xi\in\geo X$.  
It had been considered by Karpelevi\v{c} in the case of symmetric spaces,
see \cite{karpel}.

\subsubsection{Types of isometries}
\label{secparabolics}

We recall the standard classification of isometries into 
axial, elliptic and parabolic ones: 
For any isometry $\phi$ of a Hadamard space $X$ 
its displacement function $\de_{\phi}:x\mapsto d(x,\phi x)$ is convex.
$\phi$ is called semisimple if $\de_{\phi}$ attains its infimum. 
There are two types of semisimple isometries: 
$\phi$ is {\em elliptic} if the minimum is zero and has fixed points 
in this case. 
If the minimum is strictly positive then $\phi$ is {\em axial} 
and there is a non-empty family of $\phi$-invariant parallel geodesics, 
the {\em axes} of $\phi$. 
If $\de_{\phi}$ does not have a minimum then $\phi$ is called {\em parabolic}. 
The fixed point set of a parabolic isometry in $\tits X$ is non-empty 
and contained in a closed ball of radius $\pi/2$. 

\begin{dfn}
For $\xi\in\geo X$ we define the {\bf parabolic stabilizer} $P_{\xi}$
as the group consisting of all elliptic and parabolic isometries 
which preserve every horosphere centered at $\xi$.
\end{dfn}

Note that there are parabolic isometries which fix more than one point 
at infinity and do not preserve the horospheres centered at some of their
ideal fixed points. 

\begin{dfn}
\label{defunipotent}
An isometry $\phi$ of a locally compact Hadamard space $X$ is called 
{\bf purely parabolic} iff its conjugacy class accumulates at the identity.
If $Isom(X)$ is cocompact then this is
equivalent to the property that for every $\de>0$ there exist arbitrarily 
large balls on which the displacement of $\phi$ is $\leq\de$.
\end{dfn}

\subsection{Visibility Hadamard spaces}

Let $Y$ be a locally compact Hadamard space
with at least 3 ideal boundary points. 
We assume that the Tits metric on $\geo Y$ is discrete,
or equivalently, 
that $Y$ enjoys the {\em visibility property} introduced in \cite{Eberlein}:
any two points at infinity are ideal endpoints of some complete geodesic. 
Then any two distinct ideal boundary points $\xi$ and $\eta$ have Tits distance
$\pi$ 
and the family of (parallel) geodesics asymptotic to $\xi,\eta$ 
is non-empty and compact;
we denote their union by $P(\{\xi,\eta\})$. 
The visibility property is clearly inherited by closed convex subsets.
The terminology {\em visibility} is motivated by the following basic fact:

\begin{lem}
\label{visibility}
For every $y\in Y$ and every $\eps>0$ there exists $R>0$ such that the
following is true:
If $\ol{pq}$ is a geodesic segment not intersecting the ball $B_R(y)$ then
$\angle_y(p,q)\leq\eps$.
\end{lem}
\proof
See \cite{Eberlein}.
\qed

\begin{cons}
\label{compactsetofgeos}
Let $A$ be a compact subset of $\geo Y\times\geo Y\setminus Diag$.
Then the set of all geodesics $c\subset Y$ satisfying
$(c(-\infty),c(\infty))\in A$ is compact.
\end{cons}
\proof
This set $B$ of geodesic is certainly closed.
If $B$ would contain an unbounded sequence of geodesics $c_n$ then
the corresponding sequence of points $(c_n(-\infty),c_n(\infty))$ in $A$ would
accumulate at the diagonal $\Delta$,
contradicting compactness.
\qed

\begin{rem}
Visibility Hadamard spaces with cocompact isometry group
are large-scale hyperbolic in the sense of Gromov. 
\end{rem}

A sequence $(\phi_n)\subset P_{\xi}$ diverges to infinity,
$\phi_n\ra\infty$, iff $\phi_n$ converges to the constant map with value $\xi$
uniformly on compact subsets of $\geo Y\setminus\{\xi\}$.

\begin{lem}
\label{produktevonparabolischen}
Assume that for different ideal points $\xi,\eta\in\geo Y$ there are
sequences of parabolics $\phi_n\in P_{\xi}$ and $\psi_n\in P_{\eta}$ diverging
to infinity.
Then $\phi_n\psi_n$ is axial for large $n$.
\end{lem}
\proof
Let $U$ and $V$ be disjoint neighborhoods of $\xi,\eta$ respectively.
Then $\phi_n^{\pm1}(\geo Y\setminus U)\subset U$ and 
$\psi_n^{\pm1}(\geo Y\setminus V)\subset V$ 
for large $n$
which implies 
\begin{equation}
\label{parabdyn}
\al_n(\geo Y\setminus V)\subset U
\qquad\hbox{ and }\qquad
\al_n^{-1}(\geo Y\setminus U)\subset V
\end{equation}
with $\al_n=\phi_n\psi_n$. 
$\al_n$ can't be elliptic (for large $n$) because then 
$(\geo\al_n^k)_{k\in\N}$ 
would subconverge to the identity, contradicting 
(\ref{parabdyn}). 
$\al_n$ can't be parabolic either because then $(\geo\al_n^k)_{k\in\N}$ 
would converge to a constant function every where pointwise, 
which is also excluded by (\ref{parabdyn}). 
Therefore $\al_n$ is axial for large $n$.
\qed

\subsection{Buildings: Definition, vocabulary and examples}

A geometric treatment of spherical and Euclidean Tits buildings 
within the framework of Aleksandrov spaces with curvature bounded above 
has to some extent been carried through in \cite{symm}. 
We will use these results and 
for the convenience of the reader we briefly recall 
some of the basic definitions and concepts. 

\subsubsection{Spherical buildings}

A {\em spherical Coxeter complex} consists of a unit sphere $S$ 
and a finite {\em Weyl group} $W\subset Isom(S)$ 
generated by reflections at {\em walls},
i.e.\ totally geodesic subspheres of codimension 1. 
The walls divide $S$ into open convex subsets whose closures 
are the {\em chambers}. 
These are fundamental domains for the action of $W$ on $S$ 
and project isometrically to the orbit space,
the model Weyl chamber $\De_{model}=W\backslash S$. 
A {\em panel} 
is a codimension-1 face of a chamber. 

A {\em spherical building} modelled on the Coxeter complex $(S,W)$ 
is a CAT(1) space\footnote{
A {\em CAT(1) space} is a complete geodesic metric space with 
upper curvature bound 1 in the sense of Aleksandrov.}
$B$ together with an atlas of charts,
i.e.\ isometric embeddings $\iota:S\embed B$. 
The image of a chart is an {\em apartment} in $B$. 
We require that any two points are contained in an apartment 
and that the coordinate changes between charts are induced by isometries 
in $W$. 
The notions of wall, chamber, panel etc.\ transfer from the Coxeter complex
to the building. 
There is a canonical 1-Lipschitz continuous {\em accordeon map} 
$\theta_B:B\ra\De_{model}$ folding the building onto the model chamber 
so that every chamber projects isometrically. 
$\theta_B\xi$ is called the {\em type} of a point $\xi\in B$. 
$\xi$ is {\em regular} if it lies in the interior of a chamber. 

$B$ is {\em thick} if every panel is adjacent to at least 3 chanbers. 
If $B$ has no spherical de Rham factor,
i.e.\ if $W$ acts without fixed points,
then the chambers are simplices and 
$B$ carries a natural structure of a piecewise spherical 
simplicial complex. 
In this case we'll call the faces also simplices. 
A thick spherical building $B$ is called {\em irreducible}
if the corresponding linear representation of $W$ is irreducible. 
This is equivalent to the assertions
that $B$ does not decompose as a spherical join, 
and that $\De_{model}$ does not decompose. 

\bigskip
Tits originally introduced buildings 
to invert Felix Kleins Erlanger Programm 
and to provide geometric interpretations for algebraic groups,
i.e.\ to construct geometries whose automorphism groups 
are closely related to these groups. 
The simplest interesting examples of irreducible spherical buildings 
are the buildings associated to projective linear groups. 
In dimension 1, one can more generally construct a spherical building 
for every abstract projective plane, possibly with trivial group 
of projective transformations: 

\begin{ex}
Given an abstract projective plane ${\cal P}$
one constructs the corresponding 
1-dimensional irreducible spherical building $B({\cal P})$ as follows. 
There are two sorts of vertices in $B({\cal P})$: 
red vertices corresponding to points in ${\cal P}$ 
and blue vertices corresponding to lines. 
One draws an edge of length $\pi/3$ between a red and a blue vertex 
iff they are incident. 
The edges in $B({\cal P})$ correspond to lines in ${\cal P}$ 
with a marked point. 
The apartments in $B({\cal P})$, i.e.\ closed paths of length $2\pi$
and consisting of 6 edges, 
correspond to tripels of points (respectively lines) in general position. 
From the incidence properties of projective planes one easily deduces 
that any two edges are contained in an apartment
and that there are no closed paths of length $<2\pi$,
i.e.\ $B({\cal P})$ is a CAT(1) space. 
\end{ex}

Of course, 
a topological projective plane yields a topological spherical building. 

\begin{rem}[Exotic smooth projective planes]
As Bruce Kleiner pointed out 
one can produce exotic (smooth) projective planes
by perturbing a smooth projective plane, 
for instance one of the standard projective planes 
$P\R^2$, $P\C^2$ or $P\H^2$. 
\end{rem}

\subsubsection{Euclidean buildings}
\label{secEucBuil}

A {\em Euclidean Coxeter complex} consists of a Euclidean space $E$ 
and an {\em affine Weyl group} $W_{aff}\subset Isom(E)$ 
generated by reflections at {\em walls}, 
i.e.\ affine subspaces of codimension 1,
so that the image $W$ of $W_{aff}$ in $Isom(\tits E)$ is a finite reflection 
group and $(\tits E,W)$ thus a spherical Coxeter complex. 

A {\em Euclidean building} is a Hadamard space $X$ 
with the following additional structure:
There is a canonical maximal atlas of isometric embeddings 
$\iota:E\embed X$ called {\em charts} so that the coordinate changes 
are induced by isometries in $W_{aff}$. 
Any geodesic segment, ray and complete geodesic is contained in 
an {\em apartment}, i.e.\ the image of a chart. 
The charts assign to any non-degenrate segment $\ol{xy}$ 
a well-defined {\em direction} $\theta(\ol{xy})$ 
in the {\em anisotropy polyhedron} $\De_{model}$, 
the model Weyl chamber of $(\tits E,W)$. 
We request that for any two non-degenerate segments $\ol{xy}$ and $\ol{xz}$ 
the angle $\angle_x(y,z)$ takes one of the 
finitely many values which can occur in $(\tits E,W)$ 
as distance between a point of type $\theta(\ol{xy})$ and 
a point of type $\theta(\ol{xz})$.
(This is called the {\em angle rigidity property} in \cite{symm}.) 

The {\em rank} of $X$ is $dim(E)$. 
The spaces of directions $\Si_xX$ and the Tits boundary $\tits X$ 
inherite canonical spherical building structures modelled on $(\tits E,W)$. 
$X$ is {\em thick (irreducible)} if $\tits X$ is thick (irreducible). 
$X$ is called {\em discrete} 
if $W_{aff}$ is a discrete subgroup of $Isom(E)$. 
Thick locally compact Euclidean buildings are discrete 
and they carry a natural structure 
as a piecewise Euclidean simplicial complex. 

\begin{ex}
Euclidean buildings of dimension 1 are {\em metric trees},
i.e.\ spaces of infinite negative curvature in the sense that 
all geodesic triangles degenerate to tripods. 
\end{ex}

Many interesting examples 
of locally compact irreducible Euclidean buildings 
arise from simple algebraic groups 
over non-Archimedean locally compact fields with a discrete valuation.

\begin{ex}
\label{BeispA_2tilde}
Let $K$ be a locally compact field with discrete valuation, 
uniformizer $\om$, 
ring of integers ${\cal O}$ and residue field $k$. 
The Euclidean building attached to $SL(3,K)$ is constructed as follows: 
It is a simplicial complex built from isometric equilateral Euclidean 
triangles. 
The vertices are projective equivalence classes of ${\cal O}$-lattices 
in the $K$-vector space $K^3$. 
Three lattices $\La_0,\La_1,\La_2$ represent the vertices of a triangle 
if, modulo rescaling and permutation, the inclusion 
$\om\cdot\La_0\subset\La_1\subset\La_2\subset\La_0$ holds. 
$\tits X$ is isomorphic to 
the spherical building attached 
to the projective plane over $K$, 
and for any vertex $v\in X$
the space of directions $\Si_vX$ is isomorphic to 
the spherical building attached 
to the projective plane over the residue field $k$. 
\end{ex}

\begin{rem}[Unsymmetric irreducible rank-2 Euclidean buildings]
There are different locally compact fields with the same residue field,
and hence different buildings as in \ref{BeispA_2tilde} 
with isometric spaces of directions at their vertices. 
In fact one can construct uncountably many buildings such that the 
spaces of directions at their vertices are isometric to the spherical building
attached to a given projective plane. 
In this way one can obtain buildings with no non-trivial symmetry 
and their boundaries are spherical buildings attached to ``exotic''
topological projective planes. 
\end{rem}

\subsection{Locally compact topological groups}

We will make essential use of a deep result due to 
Gleason and Yamabe 
on the approximation of locally compact topological groups by Lie groups:

\begin{thm}[cf.\ {\cite[p.\ 153]{MZ}}]
\label{montgomeryzippin}
Every locally compact topological group $G$ has an open subgroup $G'$
such that $G'$ can be approximated by Lie groups
in the following sense:
Every neighborhood of the identity in $G'$ contains an invariant subgroup $H$
such that $G'/H$ is isomorphic to a Lie group.
\end{thm}

Here is a typical {\bf example} of a non-Lie 
locally compact group: 
Let $T$ be a locally finite simplicial tree and $G$ its isometry 
group equipped with the compact-open topology. 
Vertex stabilizers $Stab(v)$ 
are open compact subgroups homeomorphic to the Cantor set
and can be approximated by finite groups; 
namely every neighborhood of the identity in $Stab(v)$ 
contains the stabilizer of a finite set $V$ of vertices, $v\in V\subset T$, 
as normal subgroup of finite index. 
Other interesting examples are provided by isometry groups 
of Euclidean and hyperbolic buildings 
or more general classes of piecewise Riemannian complexes.

\section{Holonomy}
\label{holonomy}

\begin{ass}
\label{AnnWenig}
$X$ is a locally compact Hadamard space. 
$\tits X$ is a thick spherical building of dimension $r-1\geq1$. 
\end{ass}

For a unit sphere $s\subset\tits X$, 
$0\leq dim(s) <r-1$, 
we denote by $Link(s)$ 
the intersection of all closed balls $\bar B_{\pi/2}(\xi)$ centered 
at points $\xi\in s$. 
$Link(s)$ is a closed convex subset and consists of the centers 
of the unit hemispheres $h\subset\tits X$ with boundary $s$. 
Note that 
any two of these hemispheres intersect precisely in $s$ 
because $\tits X$ is a CAT(1) space. 
It won't be essential for us but is worth pointing out 
that $Link(s)$ carries a natural spherical building structure 
of dimension $dim(Link(s))=dim(\tits X)-dim(s)-1$,
compare Lemma 3.10.1 in \cite{symm}. 

For any point $\xi\in s$ we have the natural map 
\begin{equation}
\label{h1}
Link(s)\ra\Si_{\xi}\tits X
\end{equation} 
sending $\zeta$ to $\stackrel{\ra}{\xi\zeta}$. 
Both spaces $Link(s)$ and $\Si_{\xi}\tits X$ inherit a metric and 
a topology from the Tits metric and cone topology on $\tits X$,
and the injective map (\ref{h1}) is a monomorphism in the sense 
that it preserves both structures, i.e.\ it is continuous and 
a Tits isometric embedding\footnote{
Recall that the topology induced by the Tits metric is finer than the
cone topology.}. 

\begin{lem}
(\ref{h1}) maps $Link(s)$ onto $Link(\Si_{\xi}s)$. 
\end{lem}

If $dim(s)=0$ then $\Si_{\xi}s$ is empty and $Link(\Si_{\xi}s)$ 
is the full space of directions $\Si_{\xi}\tits X$. 

\proof
A direction $\stackrel{\ra}{v}\in Link(\Si_{\xi}s)$ corresponds 
to a hemisphere $h\subset\Si_{\xi}\tits X$ with boundary $\Si_{\xi}s$. 
Let $\hat\xi$ be the antipode of $\xi$ in $s$. 
Then the union of geodesics of length $\pi$ with endpoints $\xi,\hat\xi$ 
and initial directions in $h$ is a hemisphere whose center $\zeta$ 
lies in $Link(s)$ and maps to $\stackrel{\ra}{v}$. 
\qed

If $s_1,s_2\subset\tits X$ are unit spheres with
$dim(s_1)=dim(s_2)=dim(s_1\cap s_2)\geq0$ 
then for any point $\xi$ in the interior of $s_1\cap s_2$ holds 
$\Si_{\xi}s_1=\Si_{\xi}s_2=\Si_{\xi}(s_1\cap s_2)$
and the identifications 
\[ Link(s_1) \ra Link(\Si_{\xi}(s_1\cap s_2)) 
\leftarrow Link(s_2) \]
yield an {\em isomorphism}
\begin{equation}
\label{g1}
Link(s_1)  \leftrightarrow Link(s_2) 
\end{equation}
i.e.\ a cone topology homeomorphism preserving the Tits metric. 
The following lemma shows that the identification 
(\ref{g1}) does not depend on $\xi$:

\begin{lem}
For the points $\zeta_i\in Link(s_i)$ let $h_i\subset\tits X$ 
be the unit hemispheres with center $\zeta_i$ and boundary $s_i$. 
Then the points $\zeta_i$ correspond to one another under (\ref{g1}) 
iff the interiors of the hemispheres $h_i$ have non-trivial intersection.
\end{lem}
\proof
If the points $\zeta_i$ correspond to one another,
i.e.\ $\stackrel{\ra}{\xi\zeta_1}=\stackrel{\ra}{\xi\zeta_2}$, 
then the segments $\ol{\xi\zeta_i}$ initially coincide and the interiors
of the $h_i$ intersect. 
Vice versa, 
if the interiors
of the $h_i$ intersect 
then for any point $\xi$ in the interior of $s_1\cap s_2$ 
their intersection $h_1\cap h_2$ is a neighborhood of $\xi$ 
in both closed hemispheres $\bar h_i$ and therefore 
$\stackrel{\ra}{\xi\zeta_1}=\stackrel{\ra}{\xi\zeta_2}$.
\qed

\medskip
We'll now ``fill in'' the isomorphisms (\ref{g1}) 
by identifications of convex cores of cross sections of parallel sets 
in $X$. 
This will be acheived by placing different cross sections into the same 
auxiliary ambient Hadamard space, namely a space of strong asymptote classes,
so that their ideal boundaries coincide. 

Note that since $X$ has spherical building boundary, 
\ref{exuniqueminimalconvexsubset} implies that 
any apartment $a\subset\tits X$ can be filled by a $r$-flat 
$F\subset X$, i.e.\ $\geo F=a$. 
If $s\subset\tits X$ is isometric to a unit sphere 
then $s$ is contained in an apartment (by \cite[Proposition 3.9.1]{symm}) 
and hence can be filled by a flat $f\subset X$:
$\geo f=s$. 
This verifies that the parallel sets defined next are non-empty: 

\begin{dfndesc}
\label{DefParallelset}
For a unit sphere $s\subset\tits X$ we denote by $P(s)=P^X(s)$ 
the union of all flats with ideal boundary $s$. 
$P(s)$ is a non-empty convex subset and splits metrically as 
\begin{equation}
\label{ZerlegungderParallelmenge}
P(s)\cong \R^{1+\dim s}\times CS(s).
\end{equation}
The subsets $\R^{1+\dim s}\times \{point\}$ are the flats with ideal 
boundary $s$. 
$CS(s)$ is again a locally compact Hadamard space
which we call the {\bf cross section} of $P(s)$. 
For any flat $f\subset X$, 
$P(f):=P(\tits f)$ denotes its {\bf parallel set}, 
i.e.\ the union of all flats parallel to $f$,
and $CS(f):=CS(\tits f)$ denotes the cross section. 
\end{dfndesc}

Observe that $\tits CS(s)=Link(s)$.
Namely a ray in $CS(s)$ determines a flat half space in $X$ 
whose ideal boundary is a hemisphere $h$ in $\tits X$ with $\D h=s$;
vice versa, any such hemisphere in $\tits X$ can be filled by a half-flat
in $X$. 
For any point $\xi\in s$ the natural map 
$ CS(s)\ra X_{\xi} $
assigning to a point $x$ the ray $\ol{x\xi}$ 
is an isometric embedding 
because for $x_1,x_2\in CS(s)$ the triangle with vertices $x_1,x_2,\xi$ 
has right angles at the $x_i$. 

\begin{lem}
\label{EinbettungQS}
Let $s_1,s_2\subset\tits X$ be unit spheres with
$dim(s_1)=dim(s_2)=dim(s_1\cap s_2)\geq0$. 
If $\xi$ is an interior point of $s_1\cap s_2$ 
then the images of the isometric embeddings 
\begin{equation}
\label{g2}
CS(s_i)\embed X_{\xi}
\end{equation} 
have the same ideal boundary.
Furthermore the resulting identification of ideal boundaries 
coincides with the earlier identification (\ref{g1}). 
\end{lem}
\proof
Let $\zeta_i\in \tits CS(s_i)= Link(s_i)$ 
be points corresponding to each other under (\ref{g1}), 
i.e.\ $\stackrel{\ra}{\xi\zeta_1}=\stackrel{\ra}{\xi\zeta_2}$. 
The segments $\ol{\xi\zeta_i}$ initially coincide, 
i.e.\ they share a non-degenerate segment $\ol{\xi\eta}$. 
Let $r_i$ be a ray in $CS(s_i)$ asymptotic to $\zeta_i$ 
and $r'_i\subset CS(s_i)$ be the ray with same initial point 
but asymptotic to $\eta$. 
Then $r_i$ and $r'_i$ have the same image in $X_{\xi}$ under (\ref{g2})
because they lie in a flat half-plane whose boundary geodesic is asymptotic
to $\xi$. 
Since the rays $r'_1$ and $r'_2$ are asymptotic this shows that 
the images of $r_1$ and $r_2$ in $X_{\xi}$ are asymptotic rays. 
\qed

The Tits boundaries $\tits CS(s)=Link(s)$ 
contain top-dimensional unit spheres 
and \ref{exuniqueminimalconvexsubset} implies that the cross sections 
$CS(s)$ have a convex core. 

\begin{lem}
\label{keinEukFak}
If $s\subset\tits X$ is a singular sphere
then $Link(s)$ does not splitt off a spherical join factor. 
As a consequence, the convex core of $CS(s)$ has no Euclidean factor.
\end{lem}
\proof
If $Link(s)$ would have a spherical join factor then this factor would 
be contained in all maximal unit spheres in $Link(s)$. 
Hence the intersection of all apartments $a\subset\tits X$ 
with $a\supset s$ would contain a larger sphere than $s$.
This is impossible because 
$\tits X$ is a thick spherical building and 
the singular sphere $s$ is therefore an intersection of apartments. 
\qed

Fix a simplex $\tau\subset\tits X$ and choose a point $\xi$
in the interior of $\tau$. 
Then the cross sections $CS(s)$ for all singular spheres $s\supset\tau$ 
with $dim(s)=dim(\tau)$ isometrically embed into 
the same ambient Hadamard space $X_{\xi}$.
By \ref{EinbettungQS} their images have equal ideal boundaries 
and the boundary identification is given by (\ref{g1}).
According to the proof of part 2 of \ref{exuniqueminimalconvexsubset}, 
the convex cores of the $CS(s)$ are mapped to parallel layers of a flat strip 
and their boundary identifications (\ref{g1})
can be induced by isometries
which are unique in view of \ref{nontranslatable} 
and \ref{keinEukFak}. 
In this way we can compatibly identify the convex cores in consideration 
to a Hadamard space $C_{\tau}$ and there is a canonical isomorphism 
\begin{equation}
\label{d2a}
\tits C_{\tau} \buildrel\cong\over\lra \Si_{\tau}\tits X  .
\end{equation}
If $\si,\tau$ are top-dimensional simplices in the same singular sphere 
$s\subset\tits X$ 
then there is a canonical perspectivity isometry 
\begin{equation}
\label{d3}
persp_{\si\tau}:
C_{\si}\leftrightarrow C_{\tau}
:persp_{\tau\si} 
\end{equation}
because both sets are identified with the convex core of $CS(s)$.
The map of ideal boundaries induced by (\ref{d3}) 
turns via (\ref{d2a}) into an isomorphism 
\begin{equation}
\label{d4}
\Si_{\si}\tits X \leftrightarrow \Si_{\tau}\tits X
\end{equation}
(of topological buildings) which can be described inside the Tits boundary 
as follows: 
$\stackrel{\ra}{u}\in\Si_{\si}\tits X$ and 
$\stackrel{\ra}{v}\in\Si_{\tau}\tits X$
correspond to each other if they are tangent to the same 
hemisphere in $\tits X$ with boundary $s$.
(\ref{d4}) is independent of the choice of $s\supset\si\cup\tau$. 

\medskip
Let $\tau,\tilde\tau\subset\tits X$ be simplices of equal dimension 
and suppose that they are {\em projectively equivalent},
i.e.\ there exists a sequence 
$\tau=\tau_0,\dots,\tau_m=\tilde\tau$ of simplices of the same dimension 
so that any two successive simplices $\tau_i,\tau_{i+1}$ 
are top-dimensional simplices in a singular sphere. 
By composing the natural isometries (\ref{d3}), 
$C_{\tau_i}\ra C_{\tau_{i+1}}$, 
we obtain an isometry 
\begin{equation}
\label{selbst}
C_{\tau}\ra C_{\tilde\tau}
\end{equation}

\begin{dfn}
The topological space 
\[ Hol^X(\tau,\tilde\tau) \subseteq Isom(C_{\tau},C_{\tilde\tau}) \]
is defined as the closure of the subset of isometries (\ref{selbst}).
The {\bf holonomy group} 
\[ Hol(\tau)=Hol^X(\tau) \subseteq Isom(C_{\tau}) \]
at the simplex $\tau$ is defined as the topological group 
$Hol^X(\tau,\tau)$. 
\end{dfn}

For a face $\tau\subset\tits X$ 
we'd now like to relate the holonomy groupoid on the space $C_{\tau}$ 
to the holonomy groupoid on $X$.
This will be useful in the proof of \ref{bigholoindim1} 
because it allows to reduce the study of the holonomy 
action to the rank 2 case. 

Let $s,S\subset\tits X$ be unit spheres so that $s\subset S$. 
Let us denote by $s^{\perp}\subset S$ the subsphere  complementary to $s$,
i.e.\ $s^{\perp}=Link_S(s)$ and $S=s\circ s^{\perp}$. 
There are natural inclusions 
$Link(S)\subset Link(s)$ and $P(S)\subset P(s)$. 
More precisely holds
\begin{equation}
\label{i1}
Link(S)\cong Link_{Link(s)}(s^{\perp})
\end{equation}
and
\begin{equation}
\label{i2}
CS(S)\cong CS^{CS(s)}(s^{\perp})  .
\end{equation}
Assume now that the spheres $s,S$ are singular and that $\tau\subset s$ 
and $\Tau\subset S$ are top-dimensional simplices in these spheres 
so that $\tau$ is a face of $\Tau$. 
The identification $Link(s)\cong \Si_{\tau}\tits X$ 
carries $s^{\perp}$ to $\Si_{\tau}S$ and 
$Link_{Link(s)}(s^{\perp})$ to $Link(\Si_{\tau}S)$.
$core(CS(s))\cong C_{\tau}$ carries 
$core(CS(S))\cong C_{\Tau}$ to 
$core(CS^{C_{\tau}}(\Si_{\tau}S))\cong C^{C_{\tau}}_{\Si_{\tau}\Tau}$ 
and hence induces a canonical identification 
\begin{equation}
\label{i3}
C_{\Tau} \buildrel\cong\over\lra C^{C_{\tau}}_{\Si_{\tau}\Tau} . 
\end{equation}
Two faces $\Tau_1,\Tau_2\supset\tau$ are top-dimensional simplices in the same 
singular sphere $S$ iff the $\Si_{\tau}\Tau_i$ are top-dimensional simplices
in the same singular sphere in $\Si_{\tau}\tits X$. 
Let us assume that this were the case. 
Then the perspectivity $C_{\Tau_1}\leftrightarrow C_{\Tau_2}$ 
induces the perspectivity 
$C^{C_{\tau}}_{\Si_{\tau}\Tau_1} \leftrightarrow  
C^{C_{\tau}}_{\Si_{\tau}\Tau_2}$. 
We obtain an embedding 
\begin{equation}
\label{i4}
Hol^{C_{\tau}}(\Si_{\tau}\Tau) \embed Hol^X(\Tau) .
\end{equation}

\medskip
We come to the main result of this section,
namely that in the irreducible case 
the holonomy groups are non-trivial, even large: 

\begin{prop}
\label{bigholoindim1}
Suppose that, 
in addition to \ref{AnnWenig}, the spherical building $\tits X$ is
irreducible of dimension $\geq1$. 
Then for any panel $\tau\subset\tits X$
and any $\eta\in\geo C_{\tau}$, 
the parabolic stabilizer $P_{\eta}$ 
in $Hol(\tau)$ acts transitively on $\geo C_{\tau}\setminus\{\eta\}$.
\end{prop}
\proof 
Let us first consider the case $\dim\tits X=1$.
The panel $\tau$ is then a vertex $\xi$. 
The action of $Hol(\xi)$ at infinity on 
$\geo C_{\xi}\cong\Si_{\xi}\tits X$ 
can be analysed inside $\tits X$: 

\begin{sublem}
\label{sphholo2foldtrans}
For every vertex $\xi$, 
$Hol(\xi)$ acts 2-fold transitively on $\Si_{\xi}\tits X$. 
\end{sublem}
\proof
Denote by $l$ the length of Weyl chambers. 
Irreducibility implies $l\leq\pi/3$. 
Consider two vertices $\xi_1$ and $\xi_2$ of distance $2l$ 
and let $\mu$ be the midpoint of $\ol{\xi_1\xi_2}$.
Extend $\ol{\xi_1\mu\xi_2}$ 
in an arbitrary way to a (not necessarily globally minimizing)
geodesic
$\ol{\eta_1\xi_1\mu\xi_2\eta_2}$ of length $4l$. 
By irreducibility,
this geodesic is contained in an apartment $\al$
for any choice of $\eta_1$ and $\eta_2$. 
Denote by $\hat\mu$ the antipode of $\mu$ in $\al$
and let $\zeta\not\in\al$ be some neighboring vertex of $\mu$.
Then $\tangle(\zeta,\xi_i)=\pi$ and we can form the composition 
of natural maps (\ref{d4}): 
\[
\Si_{\xi_1}\tits X \ra \Si_{\zeta}\tits X \ra \Si_{\xi_2}\tits X.
\]
Varying $\eta_1,\eta_2,\zeta$ we get plenty of maps 
$\Si_{\xi_1}\tits X \ra \Si_{\xi_2}\tits X$
sending $\stackrel{\ra}{\xi_1\xi_2}=\stackrel{\ra}{\xi_1\mu}$
to $\stackrel{\ra}{\xi_2\xi_1}= \stackrel{\ra}{\xi_2\mu}$
and $\stackrel{\ra}{\xi_1\eta_1}$ to $\stackrel{\ra}{\xi_2\eta_2}$.
We can compose these and their inverses 
to obtain selfmaps of $\Si_{\xi_1}\tits X$
and see that the stabilizer of $\stackrel{\ra}{\xi_1\xi_2}$
in $Hol(\xi_1)$ acts transitively on the complement
of $\stackrel{\ra}{\xi_1\xi_2}$.
Since $\tits X$ is thick, 
$\Si_{\xi}\tits X$ contains at least three points
and it follows that $Hol(\xi)$ acts 2-fold transitively.
\qed

\medskip\no
{\em Proof of \ref{bigholoindim1} continued:}
If $P_{\eta}$ does not act transitively on $\geo C_{\tau}\setminus\{\eta\}$
then, by \ref{sphholo2foldtrans}, 
there is a non-trivial axial isometry $\al\in Hol(\xi)$ (fixing $\eta$), 
and for any $\zeta\in \geo C_{\tau}\setminus\{\eta\}$ 
there is a conjugate $\al_{\zeta}$ of $\al$ with attractive fixed point $\eta$ 
and repulsive fixed point $\zeta$. 
For $\zeta_1,\zeta_2\neq\eta$ 
the isometries $\al_{\zeta_2}^{-n}\circ\al_{\zeta_1}^n\in P_{\eta}$
subconverge to $\beta\in P_{\eta}$ 
with $\beta\zeta_1=\zeta_2$. 
This concludes the proof in the 1-dimensional case.

The general case $\dim\tits X\geq1$ can be derived:
Thanks to irreducibility, 
we can find for every panel $\tau$ an adjacent panel $\hat\tau$ 
so that $\mu:=\tau\cap\hat\tau$ has codimension 2 and
$\angle_{\mu}(\tau,\hat\tau)<\pi/2$. 
The building $\tits C_{\mu}\cong\Si_{\mu}\tits X$ 
is 1-dimensional irreducible. 
$\Si_{\mu}\tau$ is a vertex and $Hol^{C_{\mu}}(\Si_{\mu}\tau)$
acts by isometries on 
$C_{\Si_{\mu}\tau}^{C_{\mu}}\cong C_{\tau}$. 
We get an embedding 
$Hol^{C_{\mu}}(\Si_{\mu}\tau)\embed Hol^X(\tau)$
as in (\ref{i4}). 
Our result in the 1-dimensional case implies the assertion. 
\qed

\begin{ex}
If $\tits X$ is the spherical building associated to a projective plane
(with more than three points)
then $Hol(\xi)$ acts 3-fold transitive on $\Si_{\xi}\tits X$
(by ``M\"obius transformations'').
\end{ex}

\section{Rank one: Rigidity of highly symmetric visibility spaces}
\label{secrangeins}

\begin{ass}
\label{symmetriehypothese}
Let $Y$ be a locally compact Hadamard space 
with at least three ideal boundary points, 
with extendible rays, 
and which is minimal in the sense that $Y=core(Y)$. 
Suppose furthermore 
that $H\subseteq Isom(Y)$ is a closed subgroup so that
for each ideal boundary point $\xi\in\geo Y$ the parabolic
stabilizer $P_{\xi}$ in $H$ acts transitively on $\geo Y\setminus\{\xi\}$.
\end{ass}

In particular, $Y$ has the visibility property. 
For any complete geodesic $c$ we denote by $P(c)$ the {\em parallel set} of
$c$, that is, the union of all geodesics parallel to $c$.
It splits as $c\times cpt$ and contains a distinguished {\em central}
geodesic.
By \ref{produktevonparabolischen}  
there exist axial elements in $H$,
and hence the stabilizer of any central geodesic contains axial elements. 
In particular, $H$ acts cocompactly on $Y$ and $Y$ is large-scale hyperbolic
(in the sense of Gromov).
For any oriented central geodesic $c$ there is a canonical homomorphism
\begin{equation}
\label{translationsanteil}
trans:Stab(c)\ra\R
\end{equation}
given by the translational part.
Its image is non-trivial closed, so either infinite cyclic or $\R$.
The main result of this section is:

\begin{thm}
\label{rangeins}
1.\
(\ref{translationsanteil}) is surjective iff complete geodesics in $Y$ do
not branch.
In this case, $H$ is a simple Lie group, 
there exists a negatively curved symmetric space $Y_{model}$
and a homeomorphism
\[ \beta:\geo Y\buildrel\cong\over\lra \geo Y_{model} \]
which carries $H_o$ to $Isom_o(Y_{model})$: 
$\beta H\beta^{-1}=Isom_o(Y_{model})\subset Homeo(\geo Y_{model})$. 

\no
2.\
The image of (\ref{translationsanteil}) is cyclic iff $Y$ splits
metrically as $tree\times cpt$.

\no
3.\
If $T_1,T_2$ are two geodesically complete locally compact metric trees 
(with at least three ideal boundary points),
and if there are embeddings of topological groups 
$H\embed Isom(T_i)$
satisfying \ref{symmetriehypothese}, 
then there is an $H$-equivariant homothety $T_1\ra T_2$.
\end{thm}

\ref{rangeins} is a combination of the results
\ref{rangeinsindiskret},
\ref{rangeinsdiskret}
and \ref{rangeinsbaum}.

\subsection{General properties}

\begin{lem}
\label{starkasymptotischzupc}
Let $\rho:[0,\infty)\ra Y$ be a ray asymptotic to the geodesic $c$.
Then $\rho$ is strongly asymptotic to $P(c)$,
i.e.\ $d(\rho(t),P(c))\ra0$. 
\end{lem}
\proof
Assume that $\rho$ has strictly positive distance $d$ from $P(c)$.
The stabilizer of $P(c)$ contains axial elements with repulsive fixed
point $\rho(\infty)$.
Applying them to $\rho$ we can construct a geodesic at positive distance
from $P(c)$, contradicting the definition of parallel set.
\qed

\begin{lem}
Let $c$ be a geodesic and $B_{\pm}$ Busemann functions centered at the ideal
endpoints $c(\pm\infty)$.
Then the set where the 2-Lipschitz function $B_++B_-$ attains its minimum is
precisely $P(c)$.
\end{lem}
\proof
Clear.
\qed

\begin{lem}
\label{grosserwinkeldannnahe}
For every $h>0$ there exists $\al=\al(h)<\pi$ so that the following
implication holds:
If $c:\R\ra Y$ is a geodesic, $y$ a point with 
$\angle_y(c(-\infty),c(+\infty))\geq\al$
then
$d(y,P(c))\leq h$.
\end{lem}
\proof
Suppose that for some positive $h$ there is no $\al<\pi$ with this property.
Then there exist points $y_n$ of distance $\geq h$ from $P(c)$
so that $\al_n=\angle_{y_n}(c(-\infty,+\infty))\to\pi$.
(All central geodesics are equivalent modulo the action of $H$.)
This implies that there exist points $y'_n$ (on $\ol{y_n\pi_{P(c)}(y_n)}$)
so that $d(y'_n,P(c))=h$ and
$\angle_{y'_n}(c(\pm\infty),\pi_{P(c)}(y_n))\ra\pi/2$\footnote{
For a closed convex subset $C$ of a Hadamard space $X$,
$\pi_C:X\ra C$ denotes the closest point projection.}.
Since $H$ acts cocompactly,
we may assume that the $y'_n$ subconverge.
Taking a limit, we can construct a geodesic parallel to $c$
and at positive distance $h$ from $P(c)$, a contradiction.
\qed

\subsection{Butterfly construction of small axial isometries}

Consider two rays $\rho_i:[0,\infty)\ra Y$ emanating from the same point $y$
and assume that $\angle_y(\rho_1,\rho_2)<\pi$. 
Let $c_i:\R\to Y$ be extensions of the rays $\rho_i$ to complete geodesics.
We produce an isometry $\psi$ preserving the parallel set 
$P(c_1)$ by
composing four parabolic isometries:
Let $p_{i,\pm}\in P(c_i(\pm\infty))$ be the isometry which moves
$c_{i}(\mp\infty)$ to $c_{3-i}(\mp\infty)$.
Then
\[ \psi:=
p_{1,+}^{-1} p_{2,-} p_{2,+}^{-1} p_{1,-}
\]
preserves $P(c_1)$ and translates it by the displacement
\[
\de_{\psi}= 
\bigl(\sum B_{i,\pm}(y)\bigr) - \min(B_{1,+}+B_{2,-})-\min(B_{1,-}+B_{2,+})
\geq0
\]
towards $c_1(+\infty)$.
The displacement $\de_{\psi}$ is positive and $\psi$ axial iff
one of the angles $\angle_y(c_1(\pm\infty),c_2(\mp\infty))$
is smaller than $\pi$.
On the other hand, $\de_{\psi}$ is
bounded from above by twice the sum of the distances from $y$ to the parallel
sets
$Y(c_1(\pm\infty),c_2(\mp\infty))$.

\begin{lem}
\label{kleineVerschiebung}
If $\angle_y(\rho_1,\rho_2)\leq\pi-\al(h)$ then $\de_{\psi}\leq4h$.
\end{lem}
\proof
Since $\angle_y(c_1(\pm\infty),c_2(\mp\infty))\geq\al(h)$, 
\ref{grosserwinkeldannnahe} implies
$d(y,P(c_1(\pm\infty),c_2(\mp\infty)))\leq h$.
Hence 
$B_{1,\pm}(y)+B_{2,\mp}(y)-\min(B_{1,\pm}+B_{2,\mp})\leq2h$ 
and the claim follows.
\qed

\subsection{The discrete case}
\label{rankonediscretecase}

\begin{ass}
(\ref{translationsanteil}) has cyclic image:
The stabilizer in $H$ of any central geodesic has a discrete orbit
on the central geodesic.
\end{ass}

Then there is a positive lower bound for the displacement of axial
isometries in $H$. 
By \ref{kleineVerschiebung} there exists $\al_0>0$ such that:
If the rays $\rho_1$ and $\rho_2$ initiate in the same point $y$
and have angle $\angle_y(\rho_1,\rho_2)<\al_0$ then $\rho_i(\infty)$ have the
same $y$-antipodes
(i.e.\ for a third ray initiating in $y$ we have
$\angle_y(\rho,\rho_1)=\pi$ iff $\angle_y(\rho,\rho_2)=\pi$).

\begin{lem}[No small angles between rays]
\label{strahlenmitkleinemwinkel}
If the rays $\rho_1$ and $\rho_2$ initiate in the same point $y$
and have angle $<\al_0$ then
they initially coincide, i.e.\
$\rho_1(t)=\rho_2(t)$ for small positive $t$.
\end{lem}
\proof
For small positive $t$ holds
$\angle_{\rho_1(t)}(\rho_1(\infty),\rho_2(\infty))<\al_0$,
so $\rho_i(\infty)$ have the same $\rho(t)$-antipodes\footnote
{Let $x$ be a point in the Hadamard space $X$. 
Then $\xi,\eta\in\geo X$ are $x$-{\em antipodal} to each other
if there exists a geodesic passing through $x$ and asymptotic to $\xi,\eta$.}
and $\rho_2(t)=\rho_1(t)$.
\qed

\begin{lem}[Bounded Diving Time]
\label{raschesEintauchen}
If $\rho:[0,\infty)\ra Y$ is a ray asymptotic to $c$ and if
$d(\rho(0),P(c))\leq h$ then $\rho(t)\in P(c)$
for all $t\geq h/\sin(\al_0)$.
\end{lem}
\proof
We extend $\rho$ to a geodesic $c'$.
$\rho$ is strongly asymptotic to $P(c)$
(\ref{starkasymptotischzupc}).
Hence there exist $y_n\in P(c)$ tending to $\rho(\infty)$
so that the rays $\rho_n=\ol{y_nc'(-\infty)}$ Hausdorff converge to $c'$.
$\angle_{y_n}(c(-\infty),c'(-\infty))\to0$ and $\rho_n$
therefore initially lies in $P(c)$ for large $n$
(\ref{strahlenmitkleinemwinkel}). 
Outside $P(c)$ the derivative of $d(\rho_n(t),P(c))$ is $\leq-\sin(\al_0)$
whence the estimate.
\qed

\begin{cor}[Discrete Branching]
\label{verzweigung}
There exist branching complete geodesics:
Any two strongly asymptotic geodesics share a ray.
Furthermore, the set of branching points on any geodesic $c$ is discrete.
\end{cor}
\proof
The first assertion is clear from \ref{raschesEintauchen}. 
The second follows from local compactness:
Let $c_n$ be a sequence of geodesics so that
$c_n\cap c=c_n((-\infty,0])$
and the branching points $c_n(0)$ are pairwise distinct and converge.
Then, for large $n$,
the points $c_n(1)$ are uniformly separated
(by \ref{strahlenmitkleinemwinkel})
but they form a bounded subset, contradiction.
\qed

\begin{prop}[Local Conicality]
\label{lokalkegelartig}
Let $\rho:\R^+\ra Y$ be a geodesic ray,
$\si:[0,l]\ra Y$ a segment so that $\rho(0)=\si(0)$.
Then there exists $t_0>0$ so that the triangle with vertices
$\si(0),\si(t_0),\rho(\infty)$ spans a flat half-strip
and is contained in a flat strip.
\end{prop}
\proof
Denote by $\rho_t:\R^+\ra Y$ the ray emanating from $\si(t)$ and asymptotic to
$\rho$.
$\rho_t$ can be extended to a geodesic $c_t$ and there is a parallel geodesic
$c'_t$ strongly asymptotic to $\rho$.
The branch point of $c'_t$ and $\rho$ tends to $\rho(0)$ as $t\ra0$.
Discreteness of branching points on geodesics (and hence rays)
implies that $c'(t)$ passes through $\rho(0)$ for small $t$,
and $\si\restr_{[0,t]}$ lies in the flat strip bounded by $c_t$ and $c'_t$.
\qed

\begin{cons}
\label{liegtinstreifen}
Let $\rho_1,\rho_2:\R^+\ra Y$ be rays emanating from the same point $y$ and
with angle $\angle_y(\rho_1,\rho_2)=\al$.
Then $\rho_1$ can be extended to a complete geodesic $c_1$ such that
$\angle_y(\rho_2(\infty),c_1(-\infty))=\pi-\al$.
\end{cons}

\begin{cons}[Fattening half-strips]
\label{streifenerweiterung}
Let $\eta\in\geo Y$ and suppose that $\si:[0,b]\ra Y$, $0<b$, is a segment
which is contained in a complete geodesic (ray).
Assume that the ideal triangle $\De(\si(0),\si(b),\eta)$ bounds a flat
half-strip.
Then we can extend the segment $\si$ to a longer segment $\si:[a,b]\ra Y$,
$a<0$, so that the ideal triangle $\De(\si(a),\si(b),\eta)$ bounds a flat
half-strip.
\end{cons}
\proof
We assume $0<\angle_{\si(0)}(\si(b),\eta)<\pi$ 
because otherwise the claim holds trivially. 
Let $\rho:\R^+\ra Y$ be a ray extending $\si$, i.e.\
$\rho\restr_{[0,b]}\equiv\si$.
By \ref{liegtinstreifen},
we can find a geodesic $c$ extending $\rho$ and a flat strip $S$
bounded by
$c$ so that the ray $\ol{\si(0)\eta}$ is initially contained in $S$.
Then
$\angle_{\si(0)}(c(-\infty),\eta)+\angle_{\si(0)}(c(+\infty),\eta)=\pi$.
For $a<0$ sufficiently close to $0$ the ideal triangle
$\De(c(a),c(0)=\si(0),\eta)$ bounds a flat half-strip,
hence
$\angle_{c(a)}(c(b),\eta)+\angle_{c(b)}(c(a),\eta)=\pi$
and $\De(c(a),c(b),\eta)$ bounds a flat half-strip.
\qed

\begin{cor}
\label{nulloderpi}
The angle between any two rays emanating from the same point is $0$ or $\pi$.
\end{cor}
\proof
Suppose that $\rho_1,\rho_2:\R^+\ra Y$ are two rays emanating from the same
point $y$ with angle $\angle_y(\rho_1,\rho_2)=\al$.
For small $t$, the ideal triangle $\De(\rho_1(0),\rho_1(t),\rho_2(\infty))$
bounds a flat half-strip
(\ref{lokalkegelartig}).
By \ref{streifenerweiterung} and local compactness we can extend $\rho_1$ to a
complete geodesic $c_1:\R\ra Y$ so that
$\angle_{c_1(-t)}(\rho_1(\infty),\rho_2(\infty))=\al$ for all $-t\leq0$.
Since $Y$ is large-scale hyperbolic this implies that $\al=0$ or $\pi$.
\qed

\begin{prop}
\label{rangeinsdiskret}
$Y$ splits as $tree\times compact$.
\end{prop}
\proof
According to \ref{nulloderpi},
for every $y\in Y$ the union $Sun_y$ of all rays initiating in $y$ is a
minimal closed convex subset isometric to a metric tree.
\ref{exuniqueminimalconvexsubset} implies the assertion.
\qed

\begin{prop}
\label{baumdrin}
Let $Y'$ be a locally compact Hadamard space with extendible rays
and suppose that $T=core(Y')$ exists and is a metric tree. 
Then $Y'\cong T\times cpt$.
\end{prop}
\proof
The tree $T$ is locally compact and geodesically complete, so it is also
discrete.

\begin{sublem}
The nearest point projection $\pi_T:Y'\ra T$ restricts 
to an isometry on every ray $r$ in $Y'$.
\end{sublem}
\proof
We can extend $r$ to a complete geodesic $l$ and observe that 
the distance $d(\cdot,T)$ from $T$ is constant on $l$ 
because $l(\pm\infty)\in\geo T$. 
It follows that $\pi_T$ restricts on $l$ to an isometry.
\qed

\begin{sublem}
\label{winkelkannnichtwachsen}
Let $y\in Y'$ and $\xi_1,\xi_2\in\geo Y'$ so that 
$\angle_{\pi_Ty}(\xi_1,\xi_2)=\pi$. 
Then $\angle_y(\xi_1,\xi_2)=\pi$. 
\end{sublem}
\proof
For points $y_i$ on the rays $\ol{y\xi_i}$ we have 
\[d(y_1,y_2)\geq d(\pi_Ty_1,\pi_Ty_2) =
d(\pi_Ty_1,\pi_Ty) + d(\pi_Ty,\pi_Ty_2)  \]
\[= d(y_1,y)+d(y,y_2) \geq d(y_1,y_2).\]
Thus equality holds and $\angle_y(y_1,y_2)=\pi$. 
\qed

\begin{sublem}
Let $\xi\in\geo Y'$ and $c$ be a geodesic in $Y'$ not asymptotic to $\xi$. 
Then there is a point $y\in c$ with $\angle_y(l(\pm\infty),\xi)=\pi$. 
\end{sublem}
\proof
Let $y$ be the point which projects via $\pi_T$ to the center 
of the tripod in $T$ spanned by the ideal points $l(\pm\infty),\xi$
and apply \ref{winkelkannnichtwachsen}. 
\qed 

Thus any two rays in $Y'$ with same initial point 
have angle $0$ or $\pi$ 
and \ref{baumdrin} follows.
\qed

\subsubsection{Equivariant rigidity for trees}

Suppose that $T_1$ and $T_2$ are geodesically complete locally compact metric
trees with at least three boundary points,
that the locally compact topological group $H$ is embedded into their isometry
groups, $H\subseteq Isom(T_i)$,
and that the induced boundary actions of $H$ on $\geo T_i$ satisfy
\ref{symmetriehypothese}.
\begin{prop}
\label{rangeinsbaum}
Every $H$-equivariant homeomorphism $\geo T_1\ra\geo T_2$ is induced by an
$H$-equivariant homothety $T_1\ra T_2$.
\end{prop}
\proof
Maximal compact subgroups $K\subset H$ whose fixed point set on $T_i$ is a
vertex (and not the midpoint of an edge)
can be recognized from their dynamics at infinity:
There exist three ideal boundary points so that
one can map anyone to any other of them by isometries in $K$
while fixing the third. 
Adjacency of vertices can be characterised in terms of stabilizers:
The vertices $v,v'\in T_i$ are adjacent iff $Stab(v)\cap Stab(v')$
is contained in precisely two maximal compact vertex stabilizers.
It follows that there is a $H$-equivariant combinatorial isomorphism
$T_1\ra T_2$.
It is a homothety
because all edges in $T_i$ have equal length. 
\qed

\subsection{The non-discrete case}
\label{rankonenondiscretecase}

\begin{ass}
(\ref{translationsanteil}) is surjective:
The stabilizer in $H$ of any central geodesic $c$
acts transitively on $c$.
\end{ass}

\begin{lem}
\label{kleinetranslationen}
Let $G$ be an open subgroup of $H$ and $c$ a central geodesic.
Then $Stab_G(c)$ acts transitively on $c$.
\end{lem}
\proof
We choose elements $h_n\in Stab_H(c)$ with $trans(h_n)=1/n$.
They form a bounded sequence and subconverge to an elliptic element
$k\in Fix_H(c)$.
Then $(k^{-1}h_n)$ subconverges to $e$ and there exist
arbitrarily large $m\neq n$ so that
$h_m^{-1}h_n$ is axial and contained in $G$.
This shows that $Stab_G(c)$ contains axial elements with arbitrarily small
non-vanishing translational part.
\qed

\begin{cons}
\label{offeneUntergruppekokompakt}
Any open subgroup of $H$ acts cocompactly on $Y$. 
\end{cons}

\begin{prop}
\label{rangeinsindiskret}
There exist a negatively curved symmetric space $Y_{model}$,
an isomorphism $H_o\buildrel\cong\over\ra Isom_o(Y_{model})$
and an equivariant homeomorphism $\geo Y\ra\geo Y_{model}$.
\end{prop}
\proof
Suppose that
$G'\subseteq H$ is an open subgroup and
that $K$ is an invariant compact subgroup of $G'$.
\ref{kleinetranslationen}
shows that 
the $G'$-invariant non-empty closed convex subset $Fix(K)$ has full
boundary at infinity: $\geo Fix(K)=\geo Y$.
The minimality of $Y$ implies $Fix(K)=Y$ and $K=\{e\}$.
Applying \ref{montgomeryzippin}
we conclude that $H$ is a Lie group.

\begin{sublem}
$H$ has no non-trivial invariant abelian subgroup $A$.
\end{sublem}
\proof
$A$ would have a non-empty fixed point set in the geometric
compactification $\ol Y$.
If $A$ fixes points in $Y$ itself then $Fix(A)=Y$ and $A=\{e\}$
by the cocompactness of $H$ and the minimality of $Y$.
If all fixed points of $A$ lie at infinity 
then there are at most two.
This leads to a contradiction because the fixed point set of $A$ on $\geo Y$ 
is $H$-invariant, hence full or empty.
\qed

So $H$ is a semisimple Lie group with trivial center 
and $H_o\cong Isom_o(Y_{model})$ for a symmetric space
$Y_{model}$ of noncompact type and without Euclidean factor. 

\begin{sublem}
$Y_{model}$ has rank one.
\end{sublem}
\proof
If $rank(Y_{model})\geq2$ then the subgroup of translations along 
a maximal flat in $Y_{model}$
acts on $Y$ as a parabolic subgroup (because no subgroup $\cong\R^2$ in
$Isom(Y)$ can contain axial isometries) and fixes exactly one point
on $\geo Y$.
Maximal flats in $Y_{model}$ containing parallel singular geodesics yield the
same fixed point in $\geo Y$ and it follows that $H_o$ would have a fixed
point on $\geo Y$,
contradiction.
\qed

It remains to construct the equivariant homeomorphism of boundaries.
Axial isometries in $Isom(Y)$ have the property that their conjugacy class
never accumulates at the identity.
Therefore if $h\in H_o$ acts as a pure parabolic
(see definition \ref{defunipotent}) 
on $Y_{model}$ then it acts as a parabolic on $Y$.
Hence the stabilizer of $\xi_0\in\geo Y_{model}$ in $H_o$ fixes a unique point
$\xi\in\geo Y$ and we obtain an
$H_o$-equivariant, and hence continuous surjective map
$\geo Y_{model}\ra\geo Y$.
It must be injective, too, because any two stabilizers of distinct points in
$\geo Y_{model}$ generate $H_o$ but $H_o$ has no fixed point on $\geo Y$.
This concludes the proof of \ref{rangeinsindiskret}.
\qed

\begin{prop}
Complete geodesics in $Y$ don't branch.
\end{prop}
\proof
$h\in H_o$ acts as a pure parabolic on $Y$ iff it does so on $Y_{model}$.
The purely parabolic stabilizer $N_{\xi}\subset H_o$ of
$\xi\in\geo Y$
is a simply connected nilpotent Lie group and acts simply transitively on
$\geo Y\setminus\{\xi\}$.
Let $\tau\in H_o$ be any axial isometry acting on $Y$ 
with attractive fixed point $\xi$.
Then
\begin{equation}
\label{reinparabolisch}
\lim_{n\to\infty} \tau^{-n}\phi\tau^n =e .
\end{equation}
for all $\phi\in N_{\xi}$.
Let $c$ be a geodesic in $Y$ asymptotic to both fixed points of $\tau$
at infinity 
and let $\phi\in N_{\xi}$ be non-trivial.
$\tau$ acts as an isometry on the compact cross section of $P(c)$ 
and we can choose a sequence $n_k\to\infty$ so that
$d(c,\tau^{n_k} c)\to0$.
\[
d(\phi\tau^{n_k}c(0),\tau^{n_k}c(0))=d(\tau^{-n_k}\phi\tau^{n_k}c(0),c(0))\ra0
\]
implies that $\phi c$ is strongly asymptotic to $c$.
These two geodesics can't intersect because 
$\phi$ is not elliptic. 
($N_{\xi}$ has no non-trivial elliptic  elements!)
The argument shows that distinct stronlgy asymptotic geodesics are disjoint
and hence geodesics in $Y$ don't branch.
\qed

\no
{\em Proof of \ref{Rang1Analog}:}
\ref{Rang1Analog} is not much more than a reformulation of \ref{rangeins}. 
As in the proof of \ref{bigholoindim1} 
we deduce from the 2-fold transitivity of the action of $Isom(C)$ on $\geo C$ 
that the parabolic stabilizer of any $\eta\in\geo C$ 
acts transitively on $\geo C\setminus\{\eta\}$. 
Then $C$ and $Isom(C)$ satisfy assumption \ref{symmetriehypothese} 
and assertion follows from \ref{rangeins} and \ref{baumdrin}. 
\qed

\section{Geodesically complete Hadamard spaces 
with buil\-ding boun\-da\-ry}
\label{sectionhadamardspaceswithbuildingboundary}

\subsection{Basic properties of parallel sets}

\begin{ass}
\label{grundannahmeX}
$X$ is a locally compact Hadamard space with extendible rays
and $\tits X$ is a spherical building 
of dimension $r-1\geq1$. 
\end{ass}

\begin{lem}
\label{geosinflats}
Every flat half-plane $h$ in $X$ is contained in a flat plane.
\end{lem}
\proof
Let $c$ be the boundary geodesic of the flat half-plane $h$
and denote $\xi_{\pm}:=c(\pm\infty)$.
Let $\eta\in\geo h$ be so close to $\xi_+$ 
that the arc $\ol{\eta\xi_+}$ in $\tits X$ 
is contained in a closed chamber, 
and extend the ray $\ol{\eta c(0)}$ to a geodesic $c'$.
$c'$ bounds a flat half-plane $h'$ which contains $\xi_+$ in its ideal
boundary.
The canonical isometric embedding $CS(\{\xi_+,\xi_-\})\embed X_{\xi_+}$ 
sends $h$ to a ray and $h'$ to a geodesic extending this ray.
This implies that $h$ is contained in a flat plane.
\qed

\begin{cor}
\label{QuerschnittehabenerweiterbareStrahlen} 
For any flat $f\subset X$ 
the cross section $CS(f)$ is again a locally compact Hadamard space 
with extendible rays,
and $\tits CS(f)$ is a spherical building of dimension 
$\dim(\tits X)-\dim(f)$. 
\end{cor}
\proof
\ref{geosinflats} implies that for any geodesic $l$
the cross section $CS(l)$ has extendible rays. 
Now we proceed by induction on the dimension of $f$ 
using
\[ CS^{X}(f)\cong CS^{CS^{X}(f')}(CS^{f}(f')) \]
for flats $f'\subset f$.
\qed

\begin{cor}
\label{FachsliegeninApartments}
Every flat is contained in a $r$-flat. 
\end{cor}

\begin{prop}
\label{branchinggeos}
Suppose the geodesics $c_1,c_2\subset X$ have a ray $\rho$ in common.
Then there are two maximal flats whose intersection is a halfapartment.
\end{prop}
\proof
Denote $\xi:=\rho(\infty)=c_i(\infty)$ and $\xi_i:=c_i(-\infty)$.
There exist geodesics $\ga_{i}$ of length $\pi$ in $\tits X$ joining 
$\xi$ and $\xi_{i}$ so that their intersection $\ga_{1}\cap\ga_{2}$
is a non-degenerate arc $\ol{\xi\eta}$. 
The geodesics $c_{i}$ project to geodesics $\bar{c_{i}}$ in the space
of strong asymptote classes $X_{\eta}$, 
and for any $\rho(0)$-antipode $\hat\eta$ of $\eta$ the geodesics 
$\bar c_{i}$ are in fact contained in the projection to $X_{\eta}$
of the cross section $CS(\{\eta,\hat\eta\})$.  
The geodesics $\bar c_{i}$ share a ray but do not coincide 
because they have different ideal endpoints 
$\bar c_{i}(-\infty)=\stackrel{\ra}{\eta\xi_{i}}\in
\Si_{\eta}\tits X\cong\geo CS(\{\eta,\hat\eta\})$.  
We may proceed by induction on the dimension of the Tits boundary 
of the cross section until we find a flat $f$
so that $CS(f)$ 
has discrete Tits boundary and contains two geodesics whose intersection 
is a ray. 
These geodesics correspond to maximal flats in $P(f)$
with the desired property. 
\qed

\begin{reform}
\label{dannverzweigenauchdiequerschnitte}
If there are branching geodesics in $X$ then there exists a 
flat $f\subset X$ so that 
$\tits f$ is a wall in $\tits X$ and 
$CS(f)$ contains branching geodesics.
\end{reform}

\subsection{Boundary isomorphisms}
\label{secrandisos}

\begin{dfn}
Let $X'$ be another space satisfying \ref{grundannahmeX}. 
A {\bf boundary isomorphism} is a cone topology homeomorphism 
\begin{equation}
\label{randisom}
\phi:\geo X\lra\geo X'
\end{equation}
which at the same time is a Tits isometry,
i.e.\ it is an isomorphism of topological spherical buildings,
cf.\ \cite{BurnsSpatzier}. 
We denote by $Iso(\geo X,\geo X')$ the space of all 
boundary isomorphisms $\geo X\ra\geo X'$
equipped with the compact-open topology,
and by $Aut(\geo X)$ 
the topological group $Iso(\geo X,\geo X)$.
\end{dfn}

A boundary isomorphism (\ref{randisom}) 
induces for all simplices $\tau\subset\tits X$ an isomorphism
of topological buildings 
\begin{equation}
\label{differentialdesraniso}
\Si_{\tau}\tits X\lra\Si_{\phi\tau}\tits X'  .
\end{equation}
The induced homeomorphisms 
\[ 
Iso(\Si_{\tau}\tits X,\Si_{\tilde\tau}\tits X) \ra 
Iso(\Si_{\phi\tau}\tits X,\Si_{\phi\tilde\tau}\tits X) \] 
carry $Hol^X(\tau,\tilde\tau)$ to 
$Hol^{X'}(\phi\tau,\phi\tilde\tau)$ 
and thereby induce isomorphisms of topological groups
\begin{equation}
\label{induzierterholoiso}
Hol^X(\tau)\lra Hol^{X'}(\phi\tau).
\end{equation}

\begin{ass}
\label{annplusirreduzibel}
In addition to \ref{grundannahmeX} 
the building $\tits X$ is thick and irreducible. 
\end{ass}

According to \ref{QuerschnittehabenerweiterbareStrahlen}, 
$C_{\tau}$ has extendible rays.
(Extendibility of rays is inherited by subsets with full ideal boundary.) 
If $\tau\subset\tits X$ is a panel then 
by \ref{bigholoindim1}
the action of $Hol(\tau)$ on $C_{\tau}$
by isometries satisfies \ref{symmetriehypothese}
and therefore \ref{rangeins} applies. 
In the case that $\Si_{\tau}\tits X\cong\geo C_{\tau}$
is homeomorphic to a sphere, 
it can be identified with the boundary of a rank-one symmetric space 
canonically up to conformal diffeomorphism,
and $\Si_{\phi\tau}X'$ as well.
In this situation the ``differentials''
(\ref{differentialdesraniso}) 
are conformal diffeomorphisms
because they are equivariant with respect to (\ref{induzierterholoiso}). 
In the second case that $\Si_{\tau}\tits X\cong\geo C_{\tau}$
is disconnected,
$C_{\tau}$ and $C_{\phi\tau}$ are metric trees and 
(\ref{differentialdesraniso}) is conformal in the sense that it is induced 
by a homothety (\ref{rangeinsbaum}).

\medskip
The ideal boundary $\geo X$, 
equipped with the cone topology and Tits metric, 
is a compact topological spherical building. 
The cone topology can be induced by a metric 
and this allows us to apply the results from \cite{BurnsSpatzier}
on automorphism groups of topological spherical buildings. 
In particular, \cite[theorem 2.1]{BurnsSpatzier} implies: 

\begin{thm}[Burns-Spatzier]
$Aut(\geo X)$ is locally compact. 
\end{thm}

We denote by $F$ the space of chambers in $\tits X$. 
The cone topology induces a topology on $F$ which makes $F$ 
a compact space. 

\begin{lem}
There exist finitely many chambers $\si_1,\dots,\si_s$
such that the map 
\begin{equation}
\label{einbettungvonAut}
Aut(\geo X)\lra F^s\setminus Diag;
\phi\mapsto (\phi\si_1,\dots,\phi\si_s) 
\end{equation}
is proper\footnote{$Diag$ denotes the {\em generalized diagonal}
consisting of tupels with at least two equal entries.}.
\end{lem}
\proof
Choose $\si_{r+1},\dots,\si_s$ as the chambers of an apartment $a$,
and let $\tau_1,\dots,\tau_r$ be the panels of $\si_{r+1}$. 
An automorphism $\phi$ is determined by its effect on 
$a$ and the spaces $\Si_{\tau_i}\tits X$, 
because $\tits X$ is 
the convex hull of the apartment $a$ and all chambers adjacent to 
its chamber $\si_{r+1}$\footnote{
Proof:
The convex hull is a subbuilding $B'$ of maximal dimension 
\cite[prop.\ 3.10.3]{symm}. 
Since any panel is projectively equivalent to a panel $\tau_i$,
$B'$ is a neighborhood of $int(\tau)$ for any panel $\tau\subset B'$. 
We can connect an interior point of any chamber to a point in $a$ 
by a geodesic avoiding simplices of codimension $\geq2$. 
It follows that all chambers are contained in $B'$ and $B'=\tits X$.}. 
Choose for each panel $\tau_i$ a chamber $\si_i\not\subset a$
with $\si_i\cap\si_{r+1}=\tau_i$. 
Clearly (\ref{einbettungvonAut}) is continuous.
Let $(\phi_n)$ be a sequence in $Aut(\geo X)$ whose image 
under (\ref{einbettungvonAut}) is bounded,
i.e.\ does not accumulate at $Diag$. 
We have to show that $(\phi_n)$ is bounded,
respectively it suffices to show that there is a bounded subsequence.
After passing to a subsequence, 
we may assume that $\phi_n\si_i\ra\bar\si_i$ with pairwise different 
limits $\bar\si_i$. 
Denote $\bar\tau_i:=\lim\phi_n\tau_i$. 
For each $i\leq r$ the sequence of conformal homeomorphisms
$\Si_{\tau_i}\tits X\ra\Si_{\phi_n\tau_i}\tits X$ 
converges on a triple of points
(namely on $\Si_{\tau_i}(a\cup\si_i)$)
to an injective limit map and hence
subconverges uniformly to a (conformal) homeomorphism 
$\Si_{\tau_i}\tits X\ra\Si_{\bar\tau_i}\tits X$. 
It follows that $(\phi_n)$ subconverges uniformly 
to a building automorphism. 
\qed

\begin{cons}
\label{KriteriumUnbeschraenkt}
The sequence $(\phi_n)\subset Aut(\geo X)$ is unbounded 
iff there exist adjacent chambers $\si,\si'$ 
such that $\phi_n\si$ and $\phi_n\si'$ converge in $F$ 
to the same chamber.
\end{cons}

\subsection{The case of no branching}
\label{secnobranching}

A major part of the arguments in this section follows 
the lines of Gromov's proof of his Rigidity Theorem \cite{BGS} 
and the study of topological spherical buildings in \cite{BurnsSpatzier}. 

\begin{ass}
\label{annahmeunverzweigt}
$X$ is a locally compact Hadamard space with extendible rays
and $\tits X$ is a thick irreducible spherical building 
of dimension $r-1\geq1$.
Moreover 
we assume in this section that complete geodesics in $X$ do not branch.
\end{ass}

For every point $x\in X$ there is an involution
\[ \iota_x:\geo X\lra\geo X \]
which maps $\xi\in\geo X$ to the other boundary point of the unique
geodesic extending the ray $\ol{x\xi}$.

\begin{lem}
\label{enthaeltInvolutionen}
$\iota_x\in Aut(\geo X)$.
\end{lem}
\proof
The absence of branching implies that $\iota_x$ is continuous.
By \ref{dannverzweigenauchdiequerschnitte}, 
each ray emanating from $x$ is contained in a maximal flat $F$. 
$\iota_x$ restricts on the unit sphere $\tits F$ to the antipodal involution.
Hence $\iota_x$ maps every chamber isometrically to a chamber 
and is therefore 1-Lipschitz continuous with respect to the Tits distance.
The claim follows because $\iota_x^{-1}=\iota_x$. 
\qed

\ref{enthaeltInvolutionen} shows that  
the group $Aut(\geo X)$ is large. 
Our aim is to unmask it as the isometry group of a symmetric space.
Denote by $Inv$ the subgroup consisting of all products of an even number 
of involutions $\iota_{x_i}$. 

\begin{lem}
$Inv$ is path connected 
and it is contained in every open subgroup of $Aut(\geo X)$. 
\end{lem}
\proof
The map $X\times X\ra Aut(\geo X);(x_1,x_2)\mapsto\iota_{x_1}\iota_{x_2}$
is continuous 
and hence $Inv$ is path connected.
The second assertion follows in view of 
$(\iota_{x_1}\iota_{x_2})(\iota_{x_2}\iota_{x_2'})=
\iota_{x_1}\iota_{x_2'}$. 
\qed

\begin{lem}
\label{gemeinsamergegenueber}
For any two chambers $\si_1,\si_2$ in a thick spherical building $B$ 
there is a common antipodal chamber.
Refinement: 
For any two simplices of the same type\footnote{
The {\em type} of a simplex is its image under the 
canonical (accordeon) projection to the model Weyl chamber $\De_{model}$.}
there is a common antipodal simplex. 
\end{lem}
\proof
Let $\hat\si$ be a chamber antipodal to $\si_1$
and $\ga:[0,\pi]\ra B$ a unit speed geodesic avoiding codimension-2 faces
with $\ga(0)\in int(\si_2)$ and which intersects $int(\hat\si)$. 
If $\ga(\pi)\in\hat\si$ then we are done.
Otherwise let $\tau\subset\partial\hat\si$ 
be the panel where $\ga$ exits $\hat\si$. 
Since $B$ is thick, 
there exists a chamber $\hat\si'$ opposite to $\si_1$
so that $\hat\si'\cap\hat\si=\tau$. 
Let $\ga':[0,\pi]\ra B$ be a unit speed geodesic with 
$\ga'(0)=\ga(0)$, $\dot\ga'(0)=\dot\ga(0)$,
which agrees with $\ga$ up to $\hat\si$ 
and then turns through $\tau$ into the interior of $\hat\si'$. 
We repeat this procedure until it terminates after finitely steps 
and yields a chamber opposite to $\si_1$ and $\si_2$.
The refinement follows directly. 
\qed

\begin{cons}
\label{transitivauffrand}
For any simplex $\tau$, 
$Inv$ acts transitively on the compact space $F_{\tau}$
of simplices of same type as $\tau$. 
In particular, 
$Inv$ acts transitively on the compact space $F$ of Weyl chambers\footnote{
$F$ is the analog of F\"urstenberg boundary in the symmetric space case.}.
\end{cons}

Now we investigate the dynamics on $\geo X$ of elements which
correspond to translations (transvections) along geodesics
in symmetric spaces. 

\begin{lem}
\label{anzug}
Suppose that $\rho:[0,\infty)\ra X$ is a ray asymptotic to $\xi$
and that $U\subset\tits X$ is a compact set of $\xi$-antipodes. 
Then $\iota_{\rho(t)}U\ra\{\xi\}$ as $t\ra\infty$. 
\end{lem}
\proof
In every $\Si_{\eta}\tits X$, $\eta\in U$, we choose an apartment 
$\al_{\eta}$ so that the apartments 
$\geo persp_{\eta\xi}\al_{\eta}\subseteq\Si_{\xi}\tits X$ coincide.
Consider sequences $t_n\ra\infty$ and $(\eta_n)\subset U$.
We have to show that $\iota_{\rho(t_n)}\eta_n\ra\xi$. 
Let $F_n$ be a maximal flat containing the ray $\ol{\rho(t_n)\eta_n}$
and satisfying $\Si_{\eta_n}\geo F_n=\al_{\eta_n}$. 

\begin{sublem}
The family of flats $F_n$ is bounded.
\end{sublem}
\proof
Assume the contrary and, after passing to a subsequence,
that $\eta_n\ra\eta\in U$.
Denote by $a$ the unique apartment in $\tits X$ containing $\xi,\eta$
and so that $\Si_{\eta}a=\al_{\eta}$.
Let $R>0$ be large.
$F_n$ depends continuously on $t_n$ 
(by ``no branching''),
and by decreasing the $t_n$ we can acheive that 
$d(F_n,\rho(0))=R$ for almost all $n$. 
Still $t_n\ra\infty$ if $R$ is chosen sufficiently large;
namely $d(\ol{\rho(t)\eta},\rho(0))$ is bounded because there 
exists a geodesic asymptotic to $\xi$ and $\eta$. 
The $F_n$ subconverge to a maximal flat $F$ with 
$d(F,\rho(0))=R$ and $\geo F=a$.
This can't be possible for arbitrarily large $R$
because the family of flats with ideal boundary $a$ is compact,
a contradiction.
\qed

All flats arising as limits of $(F_n)$ are asymptotic to $\xi,\eta$ 
and the antipodes $\iota_{\rho(t_n)}\eta_n$ of $\eta_n$ 
in $\geo F_n$ converge to an antipode of $\eta$,
i.e.\ they converge to $\xi$.
\qed

\begin{cons}
\label{krampf}
Let $c:\R\ra X$ be a geodesic, $\xi_{\pm}:=c(\pm\infty)$
and $a_t:=\iota_{c(t)}\iota_{c(-t)}$. 
Then 
$\lim_{t\to\infty} a_t\eta=\xi_+$
iff 
$\tangle(\eta,\xi_-)=\pi$.
The convergence is uniform on compact sets of $\xi_-$-antipodes. 
\end{cons}
\proof 
By \ref{anzug}, 
$\iota_{c(-t)}\eta\ra\xi_-$ uniformly.  
Then for large $t$, $\iota_{c(-t)}\eta$ and $\xi_+$ are antipodes.
Applying \ref{anzug} again yields the claim.
\qed

Denote by $B(\xi_+,\xi_-)\subset\tits X$ 
the subbuilding defined as the union of all minimizing geodesics with 
endpoints $\xi_{\pm}$,
or equivalently, 
the union of all apartments containing $\xi_{\pm}$. 
There is a folding map (building morphism, see \cite[sec.\ 3.10]{symm}) 
$fold:\tits X\ra B(\xi_+,\xi_-)$ which is uniquely determined 
by the property that 
\[
\tangle(fold\eta,\xi_-)=\tangle(\eta,\xi_-)
\qquad\hbox{ and }\qquad
\stackrel{\lra}{\xi_-(fold\eta)}=\stackrel{\lra}{\xi_-\eta}
\]
for all $\eta\in\tits X$ with $\tangle(\eta,\xi_-)<\pi$ 
and $fold\eta=\xi_+$ if $\tangle(\eta,\xi_-)=\pi$. 
\begin{refine}
$\lim_{t\ra\infty}a_t=fold$. 
\end{refine}
\proof
By \ref{krampf} and because all $a_t$ fix the Tits neighborhood 
$B(\xi_+,\xi_-)$ of $\xi_-$ pointwise. 
\qed

\begin{prop}
\label{halbeinfach}
$Aut(\geo X)$ is a semisimple Lie group whose identity component 
has trivial center.
\end{prop}
\proof
{\em 1.\ $Aut(\geo X)$ is a Lie group:}
Let $G'\subseteq Aut(\geo X)$ be an open subgroup,
$c$ a geodesic, $\xi_{\pm}=c(\pm\infty)$ and $U_+$ a neighborhood of $\xi_+$ 
which is chosen so small that 
all points in $U_+$ with the same $\De_{model}$-direction (type) 
as $\xi_+$ are $\xi_-$-antipodes 
(using the lower semicontinuity of Tits distance). 
Suppose $H\subset G'$ is an invariant subgroup contained in the neighborhood
$\{\phi\in G':\phi\xi_+\in U_+\}$ of $e$. 
Then $H\xi_+$ consists of $\xi_-$-antipodes. 
Hence 
$H\xi_+=a_tHa_t^{-1}\xi_+=a_tH\xi_+\ra\{\xi_+\}$ 
as $t\ra\infty$,
thus $H\xi_+=\xi_+$. 
Since $Fix(H)$ is $G'$-invariant 
and convex with respect to the Tits metric 
it follows from \ref{transitivauffrand} that $Fix(H)=\geo X$ and $H=\{e\}$. 
So there are neighborhoods of the identity in $G'$ 
which don't contain non-trivial invariant subgroups. 
\ref{montgomeryzippin}
implies that $Aut(\geo X)$ is a Lie group. 

\begin{sublem}
\label{displacementpi}
Every non-trivial isometry $\phi$ of a 
thick spherical building $B$ different from a sphere
carries some point to an antipode. 
\end{sublem}
\proof
We may assume without loss of generality 
that $B$ has no spherical join factor. 
If the assertion were not true then 
$\phi$ would be homotopic to the identity
and therefore preserve every apartment and hence every simplex,
so $\phi=id$.
\qed

{\em 2.\ Semisimplicity:}
Suppose that $A$ is an invariant abelian subgroup of $Aut_o(\geo X)$.
Let $a\in A$ be a non-trivial element and 
choose a simplex $\tau_-$ such that $\tau_-$ and $a\tau_-$ are opposite 
(using \ref{displacementpi}). 
$\tau_-$ then has involution-invariant\footnote{
The type of a simplex 
is {\em involution-invariant} if its antipodal simplices have the same type,
or equivalently, if the type is fixed by the self-isometry of $\De_{model}$
which is induced by the involution of the spherical Coxeter complex.}
type. 
Let $c:\R\ra X$ be a geodesic with 
$c(-\infty)\in int(\tau_-)$ and $\tau_+$ the simplex containing $c(+\infty)$. 
Set $a_n:=\iota_{c(n)}\iota_{c(-n)}\in Inv$ and 
$b_n:=a_naa_{-n}\in A$. 
\ref{krampf} implies 
$\lim_{n\ra\infty}b_n\tau=\tau_+$
for all simplices in the open subset 
$W=\{\tau\in F_{\tau_-}:\hbox{$\tau$ and $\tau_+$ are opposite}\}$
of $F_{\tau_-}$. 
In view of \ref{gemeinsamergegenueber}, 
$W$ and the attractor $\tau_+$ are uniquely determined 
by the dynamics of $(b_n)$ and therefore are preserved by the centralizer 
of $(b_n)$ in $Aut_o(\geo X)$. 
Thus $A$ has fixed points on $F_{\tau_-}$. 
\ref{transitivauffrand} implies that 
the action of $A$ on $F_{\tau_-}$ is trivial.
The fixed point set of $A$ on $\tits X$ includes the convex hull of
all simplices in $F_{\tau_-}$ and this is the whole building $\tits X$
by irreducibility\footnote{
The convex hull of simplices of the same involution-invariant type
in the spherical Coxeter complex is a subsphere,
hence everything by irreducibility.}. 
So $A=\{e\}$. 
This shows that all abelian invariant subgroups of $Aut_o(\geo X)$ are trivial,
hence also the solvable invariant subgroups.
This finishes the proof of \ref{halbeinfach}. 
\qed

As a consequence of the proposition,
there is a symmetric space $X_{model}$ of non-compact type
and an isomorphism 
\begin{equation}
\label{gleichung}
Aut_o(\geo X)\buildrel\cong\over\lra Isom_o(X_{model}) 
\end{equation}
of Lie groups. 

\begin{lem}
\label{kompaktefixmenge}
The centralizer of every involutive boundary automorphism $\iota_x$ is
compact.
\end{lem}
\proof
Suppose that $(\phi_n)$ is an unbounded sequence in the centralizer 
of $\iota_x$.
Then there are adjacent chambers $\si,\si'$ so that 
$\lim\phi_n\si=\lim\phi_n\si'$ 
(by \ref{KriteriumUnbeschraenkt}). 
The sequence of conformal diffeomorphisms (differentials)
$\Si_{\si\cap\si'}\tits X\ra \Si_{\phi_n(\si\cap\si')}\tits X$
is unbounded and converges
everywhere except in at most one point
to a constant map 
$\Si_{\si\cap\si'}\tits X\ra \Si_{\lim\phi_n(\si\cap\si')}\tits X$. 
Denote by $s\subset\tits X$ the wall spanned by the opposite panels 
$\si\cap\si'$ and $\iota_x(\si\cap\si')$. 
It follows that for all half-apartments $h\subset\tits X$ with $\D h=s$ 
with the exception of at most one half-apartment $h_0$,
the limits $\lim\phi_n\restr_h$ exist and have the same 
half apartment $\bar h$ 
as image.
Since $\tits X$ is thick, 
we find an $\iota_x$-invariant apartment $a$ containing $s$ but not $h_0$.  
So $\phi_n\restr_a$ converges to a non-injective map $a\ra\bar h$
commuting with $\iota_x$,
i.e.\ sending antipodes to antipodes.
Such a map can't exist and we reach a contradiction.
\qed

\begin{sublem}
\label{AutosinduziertvonIsoms}
Let $X_0$ be an irreducible symmetric space. 
Every automorphism of $Isom_o(X_0)$ is the conjugation by an isometry,
i.e.\ $Isom(X_0)\cong Aut(Isom_o(X_0))$. 
\end{sublem}
\proof
$X_0=G/K$.
\qed

The involution $\iota_x\in Aut(\geo X)$ induces by conjugation 
an involutive automorphism of $Aut_o(\geo X)$, 
hence an involutive isomorphism of $Isom_o(X_{model})$ 
via (\ref{gleichung}), 
and as a consequence of \ref{kompaktefixmenge}, 
the corresponding involutive isometry of $X_{model}$ 
is the reflection at a point $\Phi(x)\in X_{model}$.
We obtain a proper continuous map
\begin{equation}
\label{AbbzuModell}
\Phi:X\lra X_{model}. 
\end{equation}
Another direct consequence is that products $\iota_x\iota_{x'}$ 
of two involutions correspond to translations (or the identity) 
in $Isom(X_{model})$. 
For any flat $F\subset X$ 
whose ideal boundary $\geo F$ is a singular sphere 
we denote by $T_F\subset Aut_o(\geo X)$ the subset 
of all $\iota_x\iota_{x'}$ with $x,x'\in F$.

\begin{lem}
\label{KorrespondenzMaxFlachs}
As a subset of $Isom(X_{model})$, $T_F$ is the group of translations 
along a flat $F^{\Phi}$ of the same dimension as $F$. 
Moreover $rank(X_{model})=r$. 
\end{lem}
\proof
Let $a\subset\tits X$ be an apartment containing $\geo F$, 
$\si$ a chamber in $a$ and $\xi_1,\dots,\xi_r$ the vertices of $\si$. 
Moreover denote by $\tau_i$ the panel of $\si$ opposite to $\xi_i$, 
and by $\hat\si,\hat\tau_i,\xi_i$ the respective antipodal objects in $a$. 
An automorphism $\phi$ of $\tits X$ which fixes $a$ pointwise is determined 
by its actions on the spaces $\Si_{\tau_i}\tits X$.
We therefore obtain an embedding
\[
Stab_{Aut(\geo X)}(a) \embed \prod_{i=1}^r Homeo(\Si_{\tau_i}\tits X) .
\]
An automorphism which fixes the subbuilding 
$\tits P(\{\xi_i,\hat\xi_i\})$ 
is determined by its action on $\Si_{\tau_i}\tits X$ alone
and we get an embedding 
\[
Stab_{Aut(\geo X)}(\tits P(\{\xi_i,\hat\xi_i\})) 
\embed Homeo(\Si_{\tau_i}\tits X).
\]
Each $\Si_{\tau_i}\tits X$ is identified with boundary of a 
rank-one symmetric space. 
$\phi\in Stab_{Aut(\geo X)}(a)$ 
acts on $\Si_{\tau_i}\tits X$ by a conformal diffeomorphism
(compare the discussion in section \ref{secrandisos}) 
which fixes at least the two point set $\Si_{\tau_i}a$. 
This diffeomorphism is hence contained in a subgroup of the conformal
group isomorphic to $\R\times cpt$.
As a consequence, 
$Stab_{Aut(\geo X)}(a)$ topologically embeds into 
a group $\cong\R^r\times cpt$
and the subgroups $H_i=Stab_{Aut(\geo X)}(\tits P(\{\xi_i,\hat\xi_i\}))$
embed into $\R\times cpt$. 
Moreover $H_i$ centralises $H_j$ for $i\neq j$. 
It follows that all translations in $Isom(X_{model})$, 
which correspond to products $\iota_x\iota_{x'}$ 
such that $x,x'$ lie on a geodesic asymptotic to $\xi_i$ and $\hat\xi_i$, 
lie in the same 1-parameter subgroup $T_i$. 
Moreover the $T_i$ commute with each other.
Since $x\mapsto\iota_x$ is proper, the first assertion follows. 

If $F$ is a maximal flat with $\geo F=a$ 
then the centralizer of $T_F$ 
is contained in $Stab_{Aut(\geo X)}(a)$
and thus contains no subgroup $\cong\R^{r+1}$.
Hence $rank(X_{model})$ can't be greater than $r$.
\qed

Consequently, 
(\ref{AbbzuModell}) sends maximal flats to maximal flats. 
Flats whose ideal boundaries are singular spheres 
arise as intersections of maximal flats 
and hence go to singular flats. 
It follows from irreducibility 
that $\Phi$ restricts to a homothety on every flat
and clearly the scale factors for restrictions to different flats agree. 
Since $X$ is geodesically complete by assumption, 
every pair of points lies in a maximal flat (\ref{FachsliegeninApartments})
and it follows that $\Phi$ is a homothety.
This concludes the proof of the main result of this section: 

\begin{thm}
\label{hauptsatzunverzweigt}
Let $X$ be a locally compact Hadamard space with extendible geodesics and
whose Tits boundary is a thick irreducible spherical building of dimension
$r-1\geq1$.
If complete geodesics in $X$ don't branch
then $X$ is a Riemannian symmetric space of rank $r$.
\end{thm}

The argument above also shows that, 
for an irreducible symmetric space $X_0$ of rank $\geq2$, 
the Lie groups $Isom(X_0)$ and $Aut(\geo X_0)$ have equal dimension 
and hence the natural embedding 
$Isom(X_0)\embed Aut(\geo X_0)$ is open 
and induces an isomorphism of identity components.
Of course, more is true:

\begin{thm}[Tits]
\label{TitsimglattenFall}
Let $X_0$ be an irreducible symmetric space of rank $\geq2$. 
Then the natural embedding
\begin{equation}
Isom(X_0)\lra Aut(\geo X_0)
\end{equation}
is an isomorphism.
\end{thm}
\proof
Let $\psi$ be an automorphism of $\geo X_0$.
We have to show that $\psi$ is induced by an isometry of $X_0$. 
$\psi$ induces an automorphism $\al$ of $Aut_o(\geo X_0)\cong Isom_o(X_0)$
which sends the stabilizer of an apartment $a$ to the stabilizer of
$\psi a$, 
i.e.\ it sends the group of translations along the flat $F_a$ 
filling in the apartment $a$ ($\geo F_a=a$) to the translations 
along $F_{\psi a}$. 
The isometry $\Psi$ inducing $\al$ (\ref{AutosinduziertvonIsoms}) 
thus satisfies $\Psi F_a=F_{\psi a}$,
i.e.\ $\geo\Psi(a)=\psi a$ for all apartments $a$ 
and it follows $\geo\Psi=\psi$.
\qed

\ref{TitsimglattenFall} implies \ref{ZusatzzuHauptsatz} in the smooth case.

\subsection{The case of branching}

\begin{ass}
\label{annahmeverzweigung}
$X$ is a locally compact Hadamard space with extendible rays
and $\tits X$ is a thick irreducible spherical building 
of dimension $r-1\geq1$.
Moreover we assume in this section that 
some complete geodesics branch in $X$.
\end{ass}

Note that now we can't expect a big group $Aut(\geo X)$
of boundary automorphisms.
There exist completely asymmetric Euclidean buildings of rank 2.
Our approach is based on the observation that nevertheless the cross sections
of all parallel sets are highly symmetric
(\ref{bigholoindim1}).

\subsubsection{Disconnectivity of F\"urstenberg boundary}
\label{eindimtitsraender}

The aim of this section is:

\begin{prop}
\label{einerdannallecantor}
If for some panel $\si$ of $B=\tits X$ the space
$\Si_{\si}B$ is totally disconnected, 
then this is true for all panels.
\end{prop}
\proof
We first consider the case when $B$ is one-dimensional.
$l$ denotes the length of a Weyl arc and  irreducibility implies $\pi/l\geq3$.
The vertices (singular points) of $B$ can be two-coloured, say blue and red,
so that adjacent vertices have different colours.
The distance of two vertices is an even multiple of $l$ iff they have the same
colour.

According to \ref{bigholoindim1}, 
the Hadamard spaces $C_{\xi}$ satisfy \ref{symmetriehypothese} 
for all vertices $\xi\in B$.
\ref{rangeins} tells that 
$\Si_{\xi}B$ is homeomorphic to a sphere of dimension $\geq1$, 
a Cantor set or a finite set with at least 3 elements 
(because $B$ is thick).  
Vertices $\xi_1,\xi_2\in B$ of the same colour are 
projectively equivalent\footnote{
For 1-dimensional spherical buildings
{\em projective equivalence} is the equivalence relation for vertices 
generated by antipodality.}
and therefore the spaces of directions $\Si_{\xi_i}B$ are homeomorphic.
If $\pi/l$ is odd then any two antipodal vertices have different colours
and the $\Si_{\xi}B$ are homeomorphic for all vertices $\xi$. 
If $\pi/l$ is even (and hence $\geq4$ by irreducibility),
we have to rule out the possibility that $\Si_{\xi}B$ is disconnected for
blue vertices $\xi$ and connected for red vertices.
Let us assume that this were the case.

\begin{sublem}
If $\Si_{\xi}B$ is a sphere for red vertices $\xi$
then $\Si_{\eta}B$ can't be finite for blue vertices $\eta$.
\end{sublem}
\proof
Assume that 
$\Si_{\xi}B$ is a sphere for red vertices $\xi$
and $\Si_{\eta}B$ is finite for blue vertices $\eta$.

{\em 1.\ Red vertices $\xi,\xi'$ of distance 4l lie in the same 
path component of the singular set $Sing(B)$:}
There exists a red vertex $\eta$ with 
$d(\xi,\eta)=d(\xi',\eta)=2l$.
$\xi,\xi',\eta$ lie in an apartment $a$
(because $4l\leq\pi$).
Let $\hat\eta$ be the antipode of $\eta$ in $a$. 
Since $\Si_{\eta}B$ is path-connected
we can continuously deform the geodesic $\ol{\eta\xi\hat\eta}$ 
to the geodesic $\ol{\eta\xi'\hat\eta}$,
so $\xi$ and $\xi'$ can be connected by a red path. 

{\em 2.\ For every red vertex $\xi_0$ 
the (red) distance sphere $S_{2l}(\xi_0)$ is path-connected:}
Let $\xi_1,\xi_2$ be red vertices with $d(\xi_i,\xi_0)=2l$.
There is a vertex $\xi'_2$ in the same path component of $S_{2l}(\xi_0)$ 
as $\xi_2$ such that $d(\xi_1,\xi'_2)=4l$.
(Deform as in 1.\ using an antipode of $\xi_0$.) 

{\em 3.\ $S_{2l}(\xi_0)$ is a manifold 
of the same dimension as $\Si_{\xi_0}B$:}
We introduce local coordinates on $S_{2l}(\xi_0)$ near $\xi$ as follows. 
Let $\eta$ be the midpoint of $\ol{\xi_0\xi}$,
i.e.\ $d(\xi_0,\eta)=d(\eta,\xi)=l$. 
Choose antipodes $\hat\xi_0$ of $\xi_0$ and $\hat\eta$ of $\eta$. 
For $\xi'\in S_{2l}(\xi_0)$ near $\xi$ the midpoint $\eta'$ 
of $\ol{\xi_0\xi'}$ is close to $\eta$ and $d(\eta',\hat\eta)=\pi$,
$d(\eta',\hat\xi_0)=d(\xi',\hat\eta)=\pi-l$. 
$\stackrel{\ra}{\hat\xi_0\eta'}$ and $\stackrel{\ra}{\hat\eta\xi'}$ 
are continuous local coordinates for $\xi'$ 
and it follows that $S_{2l}(\xi_0)$ is a manifold 
of the same dimension as $\Si_{\xi_0}B$. 

{\em 4.}
Since $\Si_{\xi_0}B$ embeds into $S_{2l}(\xi_0)$ 
it follows that $S_{2l}(\xi_0)\cong \Si_{\xi_0}B$ 
via the map $\xi\mapsto \stackrel{\ra}{\xi_0\xi}$,
and $S_{2l}(\xi_0)$ is contained in the suspension $B(\xi_0,\hat\xi_0)$. 
This implies that the cardinality of $\Si_{\eta}B$ is 2
and contradicts thickness.
\qed

For the rest of the proof of \ref{einerdannallecantor} 
we assume that $\Si_{\eta}B$ is 
a Cantor set for blue vertices $\eta$ and 
$\Si_{\xi}B$ is a sphere for red $\xi$. 

\begin{sublem}
\label{nahebei}
Let $\xi,\eta,\eta'\in B$ be distinct vertices 
(of the same color) with
$d(\xi,\eta)=d(\xi,\eta')=\pi-2l$ and
let $U$ be a neighborhood of $\xi$.
Then there exists a vertex $\xi'\in U$ satisfying
\begin{equation}
\label{sey}
d(\xi',\eta)=\pi-2l \qquad\hbox{ and }\qquad
d(\xi',\eta')=\pi.
\end{equation}
\end{sublem}
\proof
Let $\zeta$ be the vertex with $\ol{\xi\eta}\cap\ol{\xi\eta'}=\ol{\xi\zeta}$
and $\om$ the vertex on $\ol{\zeta\eta}$ adjacent to $\zeta$.
Extend $\ol{\om\xi}$ beyond $\xi$ to a geodesic $\ol{\om\hat\om}$ of length
$\pi$.
$\Si_{\om}B$ has no isolated points and we can pick a geodesic $\ga$
connecting $\om$ and $\hat\om$ so that the initial vector
$\Si_{\om}\ga$ is close to $\stackrel{\ra}{\om\zeta}$ and
the vertex $\xi'\in\ga$ with $d(\xi',\om)=d(\xi,\om)$ lies in $U$.
By construction, (\ref{sey}) holds.
\qed

\begin{sublem}
\label{inball}
Let $\ga:(-\eps,\eps)\ra Sing(B)$ be a continuous path in the 
red singular set and
$\xi$ be a red vertex so that $d(\xi,\ga(0))=\pi-2l$.
Then $d(\xi,\ga(t))=\pi-2l$ for $t$ close to $0$.
\end{sublem}
\proof
Let $\eta$ be a blue vertex adjacent to $\xi$ so that $d(\eta,\ga(0))=\pi-l$.
The set of vertices at distance $\pi-l$ from $\eta$ is open in the singular
set and so $d(\eta,\ga(t))=\pi-l$ for $t$ close to $0$.
Since $\Si_{\eta}B$ is totally disconnected we have
$\stackrel{\lra}{\eta\ga(t)}=\stackrel{\ra}{\eta\xi}$ for small $t$,
hence the claim holds.
\qed

\begin{sublem}
The path $\ga$ is constant.
\end{sublem}
\proof
It suffices to show that $\ga$ is locally constant.
There are neighborhoods $U$ of $\xi$ and $V$ of $\eta:=\ga(0)$ so that
$d(\xi',\eta')\geq\pi-2l$ for all vertices $\xi'\in U$ and $\eta'\in V$
(because the Tits distance is upper semicontinuous).
We assume without loss of generality that $\ga$ does not leave $V$.
Then \ref{inball} implies that 
$d(\xi,\ga(\cdot))\equiv\pi-2l$. 
If $\ga$ were not locally constant we could choose $t$ so that
$\eta':=\ga(t)\neq\eta$.
Applying \ref{nahebei} there exists $\xi'\in U$ so that (\ref{sey}) holds.
But \ref{inball} implies also that $d(\xi',\ga(\cdot))\equiv\pi-2l$.
Hence $d(\xi',\eta')=\pi-2l$, contradicting (\ref{sey}). 
Thus $\ga$ is locally constant. 
\qed

Hence, the set of red vertices has trivial path components. 
But since $\pi/l>2$, 
the space of directions $\Si_{\zeta}B$ for any vertex $\zeta$ continuously
embeds into the blue singular set as well as into the red singular set.
Therefore $\Si_{\zeta}B$ can't be connected for any vertex $\zeta$,
contradiction. 
Hence $\Si_{\zeta}B$ must be a Cantor set for all vertices $\zeta$.
This concludes the proof of \ref{einerdannallecantor} 
in the 1-dimensional case. 

\medskip
Without much transpiration one can deduce the assertion in the general case
$dim(B)\geq1$:
Let $\si,\tau$ be panels of the same chamber with angle
$\angle(\si,\tau)<\pi/2$.
Then the 1-dimensional topological 
spherical building $\Si_{\si\cap\tau}B$ is irreducible
and we have canonical homeomorphisms:
\[
\Si_{\si}B\cong\Si_{\Si_{\si\cap\tau}\si}\Si_{\si\cap\tau}B,
\qquad
\Si_{\tau}B\cong\Si_{\Si_{\si\cap\tau}\tau}\Si_{\si\cap\tau}B,
\]
Moreover,
$\Si_{\si\cap\tau}B$ is the ideal boundary of a Hadamard space satisfying
\ref{annahmeverzweigung} with $r=2$, 
namely of the cross section $CS(f)$ for any 
$(r-2)$-flat $f$ with $\geo f\supset\si\cap\tau$.
Therefore we can apply our assertion in the 1-dimensional case and see that
$\Si_{\si}B$ is a Cantor set if and only if $\Si_{\tau}B$ is.

Since $B$ is irreducible, for any two panels $\si,\si'$ exists a finite
sequence of panels $\si_0=\si,\si_1,\dots,\si_m=\si'$ so that any two
successive $\si_i$ are adjacent with angle less than $\pi/2$\footnote{
Otherwise we could subdivide the panels of the model Weyl chamber into two
families so that panels in different families are orthogonal; this would imply
reducibility.}.
This finishes the proof of \ref{einerdannallecantor}.
\qed

\subsubsection{The structure of parallel sets}

Consider a $(r-1)$-flat $w\subset X$ whose boundary at infinity 
is a wall in the spherical building $\tits X$. 
For any panel $\tau\subset\geo w$,
$C_{\tau}$ is canonically isometric to the convex core 
of $CS(w)$. 
By \ref{rangeins} and \ref{baumdrin}, 
the following three statements are equivalent:
\begin{itemize}
\item
$CS(w)$ is the product of a metric tree times a compact Hadamard space.
\item
$\geo CS(w)$ is homeomorphic to a Cantor set. 
\item
Some geodesics branch in $CS(w)$\footnote{
This equivalent to branching of geodesics in $C_{\tau}$.}.
\end{itemize}
By \ref{annahmeverzweigung} and \ref{dannverzweigenauchdiequerschnitte}, 
there exists a $(r-1)$-flat $w$ so that $CS(w)$ contains branching geodesics. 
$\geo w$ is a wall in $\tits X$
and for any panel $\si\subset\geo w$ we have that 
$\Si_{\si}\tits X\cong\geo CS(w)$ is a Cantor set.
\ref{einerdannallecantor} implies that $\Si_{\si}\tits X$ 
is a Cantor set for all panels $\si$ in $\tits X$ 
and hence:

\begin{lem}
$\geo CS(w)$ is $tree\times compact$ for all $(r-1)$-flats $w$. 
\end{lem}

\subsubsection{Proof of the main result \ref{hauptsatz}}
\label{BeweisvonHauptsatz}

\begin{thm}
\label{hauptsatzverzweigt}
Let $X$ be a locally compact Hadamard space with extendible rays and
whose Tits boundary is a thick irreducible spherical building of dimension
$r-1\geq1$.
If there are branching complete geodesics in $X$ then $X$ splits as the
product of a Euclidean building of rank $r$ times a compact Hadamard space.
\end{thm}
\proof
We first investigate the local structure of $X$.
For every point $x\in X$ we have the canonical 1-Lipschitz continuous 
projection
\[\theta_x:\tits X\to\Si_xX\]
which assigns to $\xi\in\tits X$ the direction
$\stackrel{\ra}{x\xi}\in\Si_xX$.
Therefore,
if $\xi,\hat\xi\in\tits X$ are $x$-antipodes,
i.e.\ $\angle_x(\xi,\hat\xi)=\pi$,
then $\theta_x$ restricts to an isometry on every geodesic 
in $\tits X$ of length $\pi$ 
connecting $\xi$ and $\hat\xi$.
By our assumption of extendible rays,
every $\xi$ has $x$-antipodes,
and it follows that $\theta_x$ restricts to an isometry on every simplex.
If $\si,\hat\si$ are open chambers in $\tits X$ which are $x$-opposite
in the sense that there exist $x$-antipodes
$\xi\in\si$ and $\hat\xi\in\hat\si$,
then $\theta_x$ restricts to an isometry on the unique apartment in $\tits X$
containing $\si,\hat\si$
and we call its image an {\em apartment} in $\Si_xX$.
If $\si_1,\si_2$ are open simplices whose $\theta_x$-images intersect
then there exists a simplex $\hat\si$ which is $x$-opposite to both $\si_i$.
It follows that the $\theta_x$-images of the spheres $span(\si_i,\hat\si)$
and therefore the $\theta_x\si_i$ coincide.
Hence the $\theta_x$-images of open simplices in $\tits X$ are disjoint
or they coincide and we call them {\em simplices} 
or {\em faces} in $\Si_xX$.

\begin{sublem}
\label{nebeneinander}
The $\theta_x$-images of adjacent chambers $\si_1,\si_2\subset\tits X$
are contained in an apartment.
(They may coincide.)
\end{sublem}
\proof
Let $\xi$ be a point in the open panel $\si_1\cap\si_2$
and $\hat\xi$ an $x$-antipode.
$\theta_x$ is isometric on the half-apartments $h_i=span(\hat\xi,\si_i)$
because it is isometric on $\si_i$ and $\angle_x(\xi,\hat\xi)=\pi$.
The union $H_i$ of the rays with initial point $x$ and ideal endpoint $\in h_i$
is a half-$r$-flat in $X$.
The $(r-1)$-flats $\D H_i$ coincide and our assumption on cross sections
of parallel sets allows two possibilities:
Either $H_1\cup H_2$ is a $r$-flat and the $\theta_x\si_i$ are adjacent
chambers in an apartment.
Or the $H_1\cap H_2$ is a non-degenerate flat strip and the $\theta_x\si_i$
coincide.
\qed

As a consequence,
the centers of adjacent chambers in $\Si_xX$ are uniformly separated and 
the compactness of $\Si_xX$ implies that the number of simplices 
in $\Si_xX$ is finite.

\begin{sublem}
Any two simplices in $\Si_xX$ are contained in an apartment.
\end{sublem}
\proof
Since in $\tits X$ any two simplices are contained in an apartment it suffices 
to prove the following statement: 
($\ast$) 
{\em If $c_1$ is a chamber contained in an apartment $a$
and $c_2$ is a chamber so that $c_2\cap a$ is a panel, 
then there exists an apartment $a'$ containing the chambers $c_1$ and $c_2$.}
The rest then follows by induction. 
To prove ($\ast$) we consider the hemisphere $h\subset a$ 
with $c_1\subset h$ and $c_2\cap a\subset\partial h$. 
Applying \ref{nebeneinander} 
to the chamber $c_2$ and the adjacent chamber in $h$ 
we see that there is a geodesic $\ga:[0,\pi]\ra\Si_xX$ 
contained in $im(\theta_x)$ which starts in $int(c_2)$,
passes through $c_2\cap a$ into $h$ and stays in $h$ for the rest of the time 
and intersects $int(c_1)$ on its way. 
The regular endpoints of $\ga$ span a unique apartment $a'$ in $\Si_xX$,
and $a'\supset c_1\cup c_2$. 
\qed

As a consequence, 
$im(\theta_x)$ is a convex, compact subset of $\Si_xX$ and hence 
itself a CAT(1)-space. 
The spherical building structure on $\tits X$ induces a 
spherical building structure on $im(\theta_x)$. 

\begin{cons}
\label{diskreteWinkel}
For any two points $\xi_1,\xi_2\in\geo X$ 
the angle $\angle_x(\xi_1,\xi_2)$ can take only finitely many values 
which depend on the types $\theta_{\tits X}\xi_i\in\De_{model}$. 
\end{cons}

Let us denote by $Sun_x$ the union of all rays emanating from $x$.

\begin{lem}
Any two sets $Sun_x$ and $Sun_y$ are disjoint or coincide.
\end{lem}
\proof
Assume that $x\neq y$ and $y\in Sun_x$. 
We pick an ideal point $\xi\in\geo X$ and
show that the ray $\ol{y\xi}$ is contained in $Sun_x$:
First we extend $\ol{yx}$ to a ray $\ol{yx\eta}$.
Then we choose a minimal geodesic connecting
$\stackrel{\ra}{yx}$ and $\stackrel{\ra}{y\xi}$ inside $im(\theta_y)$
and extend it beyond $\stackrel{\ra}{y\xi}$ to a geodesic $\al$
of length $\pi$.
Denote the endpoint by $u$.
There is a chamber $\si$ in $\Si_yX$ which contains the end of $\al$ near $u$,
and we lift $\si$ to a chamber $\tilde\si$ in $\tits X$.
$\theta_y$ restricts to an isometry on any apartment $\tilde a\subset\tits X$
which contains $\tilde\si$ and $\eta$.
Therefore $\tilde a$ bounds a flat $F$ which contains $x$ and 
a ray $\ol{y\xi'}$ with $\stackrel{\ra}{y\xi'}=\stackrel{\ra}{y\xi}$
in $\Si_yX$. 
The rays $\ol{y\xi'}$ and $\ol{y\xi}$ initially coincide 
(by \ref{diskreteWinkel}) 
and therefore 
$Sun_x\cap\ol{y\xi}$ is half-open in $\ol{y\xi}$
towards $\xi$.
Since it is clearly closed, it follows that $Sun_x$ contains the ray
$\ol{y\xi}$ and hence $Sun_y\subseteq Sun_x$.
Now the segment $\ol{xy}$ is contained in a geodesic.
I.e.\ $x\in Sun_y$ and analoguously $Sun_x\subseteq Sun_y$.
This shows that $y\in Sun_x$ iff $Sun_x=Sun_y$. 
It follows that if $z\in Sun_x\cap Sun_y$ 
then $Sun_x=Sun_z=Sun_y$,
hence the claim. 
\qed

It follows that the subsets $Sun_x$ are minimal closed convex
with full ideal boundary $\geo Sun_x=\geo X$.
Consequently they are parallel and,
by the second part of \ref{exuniqueminimalconvexsubset}, 
$X$ decomposes as a product of $Z\times compact$.
$Z$ is a geodesically complete Hadamard space 
and it remains to verify that it carries a Euclidean building structure. 
Its Tits boundary $\tits Z=\tits X$ and the spaces of directions $\Si_zZ$ 
carry spherical building structures modelled on the same Coxeter complex 
$(S,W)$
so that the maps $\theta_z:\tits Z\ra\Si_zZ$ are building morphisms,
i.e.\ they are compatible with the direction maps 
to the model Weyl chamber $\De_{model}$. 
\begin{equation}
\label{kompatibleRichtungen}
\theta_{\tits Z}=\theta_{\Si_zZ}\circ\theta_z .
\end{equation}
(The buildings $\Si_zZ$ are in general not thick.)
Choose a Euclidean $r$-space $E$, identify $\tits E\cong S$
and let $W_{aff}\subset Isom(E)$ be the full inverse image of $W$ 
under the canonical surjection $rot:Isom(E)\ra Isom(S)$. 
Up to isometries in $W_{aff}$ we can pick a canonical chart 
$E\ra F$ for every maximal flat $F\subset Z$. 
The coordinate changes will be induced by $W_{aff}$. 
Since geodesic segments are extendible they are contained in maximal flats
and in view of (\ref{kompatibleRichtungen}) 
we can assign to them well-defined $\De_{model}$-directions. 
The directions clearly satisfy the angle rigidity property 
(cf.\ section \ref{secEucBuil}) 
and we hence have a Euclidean building structure on $Z$ 
modelled on the Euclidean Coxeter complex $(E,W_{aff})$. 
(If one wishes, one can reduce the affine Weyl group and obtain a 
canonical thick Euclidean building structure.) 
This concludes the proof of \ref{hauptsatzverzweigt}. 
\qed

\medskip\no
{\em Proof of \ref{hauptsatz}:}
Put in \ref{hauptsatzunverzweigt} and \ref{hauptsatzverzweigt}.
Stir gently.
\qed

\subsection{Inducing boundary isomorphisms by homotheties: 
Pro\-of of \ref{ZusatzzuHauptsatz}}
\label{fliege}

\no
{\em Proof of \ref{ZusatzzuHauptsatz}:}
By \ref{hauptsatz},
$X$ and $X'$ are symmetric spaces or Euclidean buildings.
The F\"urstenberg boundary of $X$ is a Cantor set
iff $X$ is a Euclidean building. 
Hence $X,X'$ are either both symmetric or both buildings.
The assertion in the symmetric case is the content of 
\ref{TitsimglattenFall}. 

We may therefore assume that $X$ and $X'$ are 
thick irreducible Euclidean buildings of rank $r\geq2$. 
Then for any flat $f$ the cross section $CS(f)$ is a Euclidean building 
of rank $r-\dim(f)$,
and has no Euclidean factor if $f$ is singular. 
For all geodesics $l$
the canonical embeddings
$CS(l)\embed X_{l(\pm\infty)}$ 
of cross sections into spaces of strong asymptote classes are now 
surjective isometries, 
and for every $\xi\in\geo X$, 
$C_{\xi}\cong X_{\xi}$ is a Euclidean building 
of rank $r-1$ which splitts off a Euclidean de Rham factor 
of dimension $\dim(\tau_{\xi})$ where $\tau_{\xi}$ denotes 
the simplex containing $\xi$ as interior point. 
In particular, 
for all panels $\tau\subset\tits X$, $C_{\tau}$ is a metric tree.  
As explained in section \ref{secrandisos},
the differentials (\ref{differentialdesraniso}) of (\ref{RandIso}) 
are boundary maps of homotheties 
\begin{equation}
C_{\tau}\lra C_{\phi\tau} 
\end{equation}
and these commute with the system of natural perspectivity identifications 
(\ref{d3}). 
The assertion of \ref{ZusatzzuHauptsatz} 
follows if we can pin down every vertex of $X$ by data at infinity. 
This is acheived by the following {\bf bowtie construction} 
suggested by Bruce Kleiner: 
A {\em bowtie} $\bowtie$ 
consists of a pair of opposite chambers $\si_{\bowtie}$ and 
$\hat\si_{\bowtie}$, 
of vertices $y_i\in C_{\tau_i}$ for each panel $\tau_i\subset\si_{\bowtie}$ 
and vertices $\hat y_i\in C_{\hat\tau_i}$ for the opposite panels $\hat\tau_i$ 
so that 
$persp_{\tau_i\hat\tau_i}y_i=\hat y_i$ 
holds. 
$\bowtie$ determines a vertex in $X$ as follows: 
$\si_{\bowtie}$ and $\hat\si_{\bowtie}$ 
are contained in the ideal boundary of a 
unique maximal flat $F_{\bowtie}\subset X$
and every pair $y_i,\hat y_i$ determines a wall $w_i\subset F_{\bowtie}$. 
The $r$ walls $w_i$ intersect in a unique vertex $x_{\bowtie}$. 
We say loosely that $\bowtie$ is {\em contained} in the flat $F_{\bowtie}$. 
We call two bowties $\bowtie$ and $\bowtie'$ pre-adjacent if
$\si_{\bowtie}\cap\si'_{\bowtie}=\tau_r$, 
$\hat\si_{\bowtie}=\hat\si'_{\bowtie}$
and $\hat y_i=\hat y_i'$ for all $i$. 
(Then also $y_r=y_r'$ holds.) 
There is an obvious involution on the space of bowties 
and an equally obvious action of the permutation group $S_r$
and we call two bowties {\em adjacent} 
if they are pre-adjacent modulo these operations.
Adjacent bowties determine the same vertex. 
Adjacency spans an equivalence relation on the set of bowties
which we denote by ``$\sim$''.

\begin{lem}
\label{fliegeecke}
$\bowtie\sim\bowtie'$ iff $x_{\bowtie}=x_{\bowtie'}$.
\end{lem}
\proof
Clearly 
$\bowtie\sim\bowtie'$ implies $x_{\bowtie}=x_{\bowtie'}$.
To prove the converse, let us assume that $x_{\bowtie}=x_{\bowtie'}$. 
We start with a special case: 

\begin{sublem}
\label{FliegeninWohnung}
If $\bowtie$ and $\bowtie'$ lie in the same apartment then 
$\bowtie\sim\bowtie'$. 
\end{sublem}
\proof
It is enough to check the case when $\si_{\bowtie}$ 
and $\si_{\bowtie'}$ share a panel, 
i.e.\ without loss of generality 
$\tau_1=\tau'_1$, $\hat\tau_1=\hat\tau'_1$,
$y_1=y_1'$ and $\hat y_1=\hat y_1'$. 
Since $y_1,\hat y_1$ are vertices 
there exists a half-$r$-flat $H\subset X$ so that 
$H\cap F=\D H=w_1$. 
If $\bowtie''$ is adjacent to $\bowtie$ and 
$\bowtie'''$ is adjacent to $\bowtie'$ 
so that 
$\si_{\bowtie''}=\si_{\bowtie'''}\subset\geo H$ 
then $\bowtie''$ and $\bowtie'''$ are adjacent.
So $\bowtie\sim\bowtie'$. 
\qed

\begin{sublem}
\label{ersetzeFliege0}
Let $\hat\si$ be a chamber in $\tits X$. 
Then there exists a bowtie $\bowtie''\sim\bowtie$ 
so that $\hat\si_{\bowtie''}=\hat\si$.
\end{sublem}
\proof
It is enough to treat the case when $\hat\si$ is adjacent to 
$\hat\si_{\bowtie}$ in $\tits X$. 
After replacing $\bowtie$ by an equivalent bowtie 
(e.g.\ contained in the same maximal flat) 
we may assume that 
$\theta_{x_{\bowtie}}\hat\si$ and
$\theta_{x_{\bowtie}}\si_{\bowtie}$ 
are opposite chambers in $\Si_{x_{\bowtie}}X$. 
Then we can choose $\bowtie''$ adjacent to $\bowtie$.
\qed

We refine the previous sublemma: 

\begin{sublem}
\label{ersetzeFliege}
Let $F$ be a maximal flat and $\bowtie$ a bowtie with 
$\si_{\bowtie}\subset\geo F$ and $x_{\bowtie}\in F$. 
Let $\hat\si$ be any chamber in $\tits X$. 
Then there exists another bowtie $\bowtie''\sim\bowtie$ 
so that $\hat\si_{\bowtie''}=\hat\si$ and $\si_{\bowtie''}\subset\geo F$. 
\end{sublem}
\proof
Again, 
we may assume without loss of generality that the chamber 
$\hat\si$ is adjacent to $\hat\si_{\bowtie}$. 
If $\theta_{x_{\bowtie}}\hat\si$ is opposite to the chamber 
$\theta_{x_{\bowtie}}\si_{\bowtie}$ in $\Si_{x_{\bowtie}}X$ 
then we can choose $\bowtie''$ adjacent to $\bowtie$. 
Otherwise 
let $\si\subset\geo F$ be the chamber adjacent to $\si_{\bowtie}$ so that 
$\theta_{x_{\bowtie}}\si$ is opposite to
$\theta_{x_{\bowtie}}\hat\si$
and denote by $\bowtie''$ the bowtie with 
$\si_{\bowtie''}=\si$, $\hat\si_{\bowtie''}=\hat\si$ 
and $x_{\bowtie''}=x_{\bowtie}$. 
Then $\bowtie''$ is equivalent to a bowtie contained in $F_{\bowtie}$
and hence to $\bowtie$.
\qed

To finish the proof of \ref{fliegeecke}
we can first replace $\bowtie'$ by an equivalent bowtie
so that $\si_{\bowtie}=\si_{\bowtie'}$ 
(\ref{ersetzeFliege0})
and then replace it in a second step 
so that $\bowtie'$ and $\bowtie$ lie in the same apartment
(\ref{ersetzeFliege}). 
Hence $\bowtie'$ and $\bowtie$ are equivalent
(\ref{FliegeninWohnung}). 
\qed

It follows that 
equivalence classes of bowties in $X$ correspond to vertices. 
Since (\ref{RandIso}) induces a map between the spaces of bowties in $X$ 
and $X'$ which preserves the equivalence relation ``$\sim$'',
it thereby induces a map $\Phi:Vert(X)\ra Vert(X')$ on vertices. 
$\Phi$ maps all vertices in a singular flat $f\subset X$ 
to the vertices of a singular flat $f^{\Phi}\subset X'$ 
so that $\phi(\geo f)=\geo f^{\Phi}$. 
Since $X$ and $X'$ are irreducible buildings,
$\Phi$ extends to a homothety $\Phi:X\ra X'$
and $\geo\Phi=\phi$.
This concludes the proof of \ref{ZusatzzuHauptsatz}. 
\qed

\subsection{Extension of Mostow and Prasad Rigidity to singular spaces 
of nonpositive curvature: 
Proof of \ref{verallgMostow}} 

\no
{\em Proof of \ref{verallgMostow}:}
We argue as Mostow \cite{Mostow}.
A $\Ga$-{\em periodic} flat is a maximal flat whose stabilizer in $\Ga$
acts cocompactly.
Due to results of Borel and Ballmann-Brin,
$\Ga$-periodic flats lie dense in the space of all flats in $X_{model}$.
By our assumption,
there is a $\Ga$-equivariant continuous map
\[ \Phi:X_{model}\lra X .\]
It is a quasi-isometry and carries $\Ga$-periodic flats in $X_{model}$ to
$\Ga$-periodic quasi-flats in $X$ with uniform quasi-isometry constants.
If a quasi-flat is Hausdorff close to a flat then it lies in a tubular
neighborhood of this flat whose radius is uniformly bounded in terms of the
quasi-isometry constants
(\cite[Lemma 13.2]{Mostow} for symmetric spaces and \cite{symm} for
buildings).
Density and uniformity imply that $\Phi$ maps every flat in $X_{model}$ 
uniformly close to a flat in $X$ 
and with this information 
one can construct a $\Ga$-equivariant boundary isomorphism
\[ \Phi_{\infty} :\geo X_{model}\lra \geo X .\]
By \ref{hauptsatz} $X$ is a symmetric space or Euclidean building, 
and by \ref{ZusatzzuHauptsatz},  
after suitably rescaling the irreducible factors of $X_{model}$,
$\Phi_{\infty}$ is induced by a $\Ga$-equivariant isometry
$X_{model}\ra X$.
\qed


\begin{thebibliography}{BaBE}
\addcontentsline{toc}{section}{Bibliography}

\bibitem[Al57]{Aleks}
A.D.\ Aleksandrov, 
{\em \"Uber eine Verallgemeinerung der Riemannschen Geometrie,} 
Schriften des Forschungsinstituts f.\ Mathematik 1 (1957), 33-84. 

\bibitem[Ba85]{Rangstarrheit}
W.\ Ballmann,
{\em Nonpositively curved manifolds of higher rank,}
Ann.\ Math.\ 122 (1985), 597-609.

\bibitem[Ba95]{Ballmann}
W.\ Ballmann,
{\em Lectures on spaces of nonpositive curvature},
DMV-Seminar notes, vol.\ 25, Birkh\"auser 1995.

\bibitem[BGS85]{BGS}
W.\ Ballmann, M.\ Gromov, V.\ Schroeder,
{\em Manifolds of Nonpositive Curvature,}
Birkh\"auser 1985.

\bibitem[BS87]{BurnsSpatzier}
K.\ Burns and R.\ Spatzier,
{\em On topological Tits buildings and their classification,}
Publ.\ IHES 65 (1987), 35-59.

\bibitem[EO73]{Eberlein}
P.\ Eberlein and B.\ O'Neill,
{\em Visibility manifolds,}
Pacific J.\ Math.\ 46 (1973), 45-109. 

\bibitem[Eb88]{symmdiff}
P.\ Eberlein, 
{\em Symmetry Diffeomorphism Group of a Manifold of Nonpositive Curvature II,} 
Indiana U.\ Math.\ J.\ 37, no.\ 4 (1988), 735 -752. 

\bibitem[Gl52]{Gleason}
A.\ Gleason,
{\em Groups without small subgroups,}
Ann.\ Math.\ 56 (1952), 193-212.

\bibitem[Gr87]{HypGrp}
M.\ Gromov,
{\em Hyperbolic Groups,}
in: Essays in Group Theory,
MSRI Publications vol.\ 8, Springer 1987.

\bibitem[Gr93]{AsyInv}
M.\ Gromov,
{\em Asymptotic invariants for infinite groups,}
in: Geometric group theory, 
London Math.\ Soc.\ lecture note series 182, 1993.

\bibitem[Ka]{karpel}
F.\ I.\ Karpelevi\v{c},
{\em The geometry of geodesics and the eigenfunctions of the Beltrami-Laplace
operator on symmetric spaces,}
Amer.\ Math.\ Soc.\ Translations 14 (1985), 51-199.

\bibitem[KL96]{symm}
B.\ Kleiner and B.\ Leeb,
{\em Rigidity of quasi-isometries for symmetric spaces and Euclidean
buildings,}
preprint, Bonn april 1, 1996,
to appear in Publications IHES.

\bibitem[Le95]{thesis}
B.\ Leeb,
{\em 3-manifolds with(out) metrics of nonpositive curvature,} 
Inventiones math.\ 122 (1995), 277-289. 

\bibitem[Ma91]{Margulis}
G.\ Margulis,
{\em Discrete subgroups of semisimple Lie groups,}
Springer 1991.

\bibitem[MZ55]{MZ}
D.\ Montgomery and L.\ Zippin,
{\em Topological transformation groups,}
Interscience Publishers, New York 1955.

\bibitem[Mos73]{Mostow}
G.\ D.\ Mostow,
{\em Strong rigidity of locally symmetric spaces,}
Princeton UP 1973.

\bibitem[Pra79]{Prasad}
G.\ Prasad,
{\em Lattices in semisimple groups over local fields,}
Studies in algebra and number theory, 
Adv.\ Math.\ Suppl.\ Studies 6 (1979). 

\bibitem[Ron89]{Ronan}
M.\ Ronan, 
{\em Lectures on Buildings,} 
Perspectives in Mathematics vol.\ 7,
Academic Press 1989.

\bibitem[Ti74]{Tits}
J.\ Tits,
{\em Buildings of Spherical Type and Finite BN-Pairs,}
LNM 386, Springer 1974 (2nd ed.\ 1986). 

\bibitem[Ya53]{Yamabe}
H.\ Yamabe,
{\em A generalization of a theorem of Gleason,}
Ann.\ Math.\ 58 (1953), 351-365.

\end{thebibliography}
\end{document}